\documentclass[a4paper,10pt]{article}
\usepackage[utf8]{inputenc}
\usepackage{textcomp}
\usepackage{graphicx}
\usepackage{graphics, amsmath,amsthm}
\usepackage{bbm}
\usepackage{MnSymbol}
\DeclareMathAlphabet\mathbb{U}{msb}{m}{n}
\usepackage[T1]{fontenc}
\usepackage{float}
\usepackage{dsfont}
\newcommand{\vep}{\varepsilon}
\renewcommand{\P}{\mathbb{P}}
\newcommand{\N}{\mathbb{N}}
\newcommand{\Z}{\mathbb{Z}}
\newcommand{\R}{\mathbb{R}}
\newcommand{\E}{\mathbb{E}}
\newcommand{\ra}{\rightarrow}
\newcommand{\mc}[1]{\mathcal{#1}}
\newcommand\1{\mathrm{1}}
\newcommand{\indic}[1]{\1_{\{#1\}}}

\newtheorem{Lemma}{Lemma}[section]

\newtheorem{Theorem}[Lemma]{Theorem}
\newtheorem{Conjecture}[Lemma]{Conjecture}

\newcommand{\Qed}{\qed \medskip}
\newcommand{\iB}{\mathring{B_{L}}}
\newcommand{\io}{\mathring{\omega}}
\newcommand{\iO}{\mathring{\Omega}}

\newcommand{\del}{\partial}

\newcommand{\about}{\approx}

\title{The voter model chordal interface in two dimensions}
\author{Mark Holmes$^{a}$, Yevhen Mohylevskyy$^{a}$, Charles M. Newman$^{b}$\\
\begin{scriptsize}                                                                        
$^{a}$  University of Auckland, Auckland, New Zealand 
\end{scriptsize}
\\
\begin{scriptsize} 
$^{b}$ Courant Institute of Mathematical Sciences and NYU-Shanghai
\end{scriptsize}
}

\begin{document}

\maketitle

\begin{abstract}
Consider the voter model on a box of side length $L$ (in the triangular lattice) with boundary votes fixed forever as type 0 or type 1 on two different halves of the boundary. Motivated by analogous questions in percolation, we study several geometric objects at stationarity, as $L\rightarrow \infty$. One is the interface between the (large -- i.e., boundary connected) 0-cluster and 1-cluster. Another is the set of large ``coalescing classes'' determined by the coalescing walk process dual to the voter model. 
\end{abstract}
\nocite{}
\section{Introduction}
\label{sec:intro}
In this section we motivate our study of the (two-dimensional) voter model and its dual coalescing walks through their connection with a number of percolation models. In Section 2, we report on numerical results for the dimension of a natural ``chordal interface'' of the voter model.  In Section 3 we give rigorous (and a few numerical) results on the large coalescing classes for coalescing walks (where vertices $x$ and $y$ in a box are in the same class if their walks coalesce before hitting the boundary). In the appendix, more details about our numerical results are provided.

\begin{figure}
\begin{center}
\hspace{1cm}
\includegraphics[scale=.5]{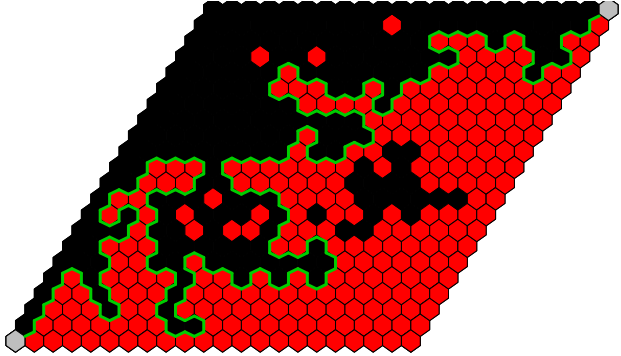}
\caption{A possible percolation configuration on $B_L$ with $L=22$, with the exploration path shown in green.}
\label{fig:small_perc}
\end{center}
\end{figure}

Among the most important breakthroughs in statistical physics and probability in the last two decades is the work by Schramm and coauthors \cite{LawSchWer04,Sch00,SchShe05} and Smirnov \cite{Sm01, Sm01b, Sm06} identifying (or conjecturing) members of the Schramm-Loewner Evolution family of random curves as the scaling limits of various random walks and interfaces in two-dimensional spin systems.  In particular Smirnov \cite{Sm01,Sm01b} (see also Camia and Newman's paper \cite{CN07}) has shown that the scaling limit of critical site-percolation on the triangular lattice $\mathbb{T}$ is SLE$_6$.  To give a rough description of one version of this statement, take a rhombic box $B_{L}$ (containing $L \times L$ vertices) in the triangular lattice in two dimensions.  Label the sides clockwise starting from the southwest corner as $\partial_1,\partial_2,\partial_3,\partial_4$.
A percolation configuration on $B_{L}$ is an element $\omega=(\omega_x)_{x\in B_L}$ of $\Omega_L=\{0,1\}^{B_L}$ defined as follows.  Fix the vertices in $\partial_1$ and $\partial_2$ to  have value 0 (or black or closed) and those in $\partial_3$ and $\partial_4$ to be 1 (or red or open).  In the interior $\iB$, set each vertex to be (independently) 0 or 1, with probability $1/2$ each -- see Figure \ref{fig:small_perc}.
There exists a unique simple path $Z_L(\omega)$ of length $|Z_L|$ from the southwest corner following edges in the dual hexagonal lattice to the opposite corner that keeps black/closed vertices on the left and red/open vertices on the right. $Z_L$ is often referred to as the exploration path; we will also call it the chordal interface. As $L\nearrow \infty$ the law of $Z_L$, after rescaling, converges weakly to a probability measure on continuous paths that is the law of {\em chordal SLE$_6$} \cite{CN07,Sm01} in a rhombic domain.  One can use this to prove (see \cite[Prop. 2]{NolWer09}) that $\E[|Z_L|]\approx L^{7/4}$, where
\begin{align}
\label{about}
f(L)\about g(L) \qquad \iff \qquad \frac{\log(f(L))}{\log(g(L))}\ra 1, \quad \text{as }L \ra \infty.
\end{align}
In general for $\beta>0$ we have that $f(L)\approx L^\beta$  if and only if 
\begin{align}
f(L)=L^\beta\ell(L), \quad \text{where} \quad \frac{\log \ell(L)}{\log L}\rightarrow 0, \quad \text{as } L \rightarrow \infty.\label{approxbeta} 
\end{align}
Note that $\omega\in \Omega_L$ uniquely determines the path $Z_L(\omega)$.  One can therefore ask about the limiting behavior of $Z_L$ when the configurations $\omega$ are generated by some other process (i.e., not i.i.d.~critical site percolation) in the interior $\iB$.  In the case of the Ising model (where the states at two sites are not independent) at the critical temperature, Smirnov \cite{Sm06} has identified  that the limiting probability measure is instead chordal {\em SLE$_3$}.
We are interested in the limiting behavior of $Z_L$ when the law of the configuration $\omega$ is the stationary distribution of the voter model (or related models) on $\iB$.  

\subsection{The voter model}
\label{sec:voter}
%Take a rhombic box $B_{L}$ of size $L \times L$ in the triangular (2 dimensional) lattice. On this graph we run a stochastic process $V_{t}$ called the voter model, with fixed boundary conditions.  Each vertex $v \in B_{L}$ takes one of the two possible values (colors or opinions) 0 or 1, so the state space for the process $V_{t}$ is $\{ 0,1 \} ^{B_{L}}$. Colors for vertices $v \in \partial B_{L}$, the boundary of $B_{L}$, are fixed to be 0 for two adjacent sides with 120\textdegree angle in between, 1 for the other two, and do not change in time. Colors for $v \in \mathring{B_{L}}$, the interior of $B_{L}$, are chosen according to some initial distribution on  $\{0,1\}^{\ \mathring{B_{L}}}$. 

%The process $\{V_{t}\}_{t\ge 0}=\{V_{t}\}_{t\ge 0}(L)$ evolves in continuous time in the following way. 
In this section we define our primary model of interest, on $B_{L}$ as described above, with  boundary states set as 0 on one pair of adjacent sides and 1 on the other pair, while the law of the interior states $\io=(\omega_x)_{x \in \iB}$ is the stationary measure for the voter model $\{V_t\}_{t\ge 0}$ on $\iB$, as follows.

Each  $v \in \mathring{B_{L}}$ has its own independent Poisson clock (a Poisson process $\Gamma_v$) of rate 1. When the clock of a vertex $v$ rings we update the state $V_t(v)$ of $v$ by choosing one of its six neighbors uniformly at random and adopting the state of the chosen neighbor.  Note that the neighbor may be one of the vertices in the boundary $\partial B_{L}=B_{L}\setminus \mathring{B_{L}}$ whose state is fixed. 
Defined this way, $V_{t}$ is an irreducible Markov process with finite state space $\iO=\{0,1\}^{\iB}$, and therefore it has a unique invariant distribution.  We will write $V_{\infty}^L$ for a random configuration sampled from this invariant distribution.

%\subsection{Coalescing random walks and duality}

The process admits a well-known graphical representation (due to T.E. Harris \cite{Harris78}) which we now review. For each $v \in B_{L}$, we draw a positive half line (representing time) in the third dimension, and on it we mark the times of Poisson clock rings of that vertex. Each mark on a time line represents a state update event which also has an arrow from $v$ to the uniformly chosen neighbor whose state is adopted. The lines of the boundary vertices have arrow marks to them, but not from them, as those states are fixed. 

Fix an initial configuration $V_0=\{V_0(x)\}_{x\in \iB}$.  To determine the state of a vertex $v\in \iB$ at time $t$ we start at height/time $t$ on the time line corresponding to $v$ and follow it down until we reach height/time 0 or we encounter an outgoing arrow (whichever comes first) at height $t'\in (0,t)$. If we meet an outgoing arrow we follow it to the time line of a neighboring vertex $v'$ and proceed as before, following this time line down from height $t'$ until reaching height/time 0 or an outgoing arrow.  We stop this procedure when we reach a boundary vertex or height 0 on some time line.  Thus from any $v\in \iB$ and $t>0$ the path followed corresponds to a continuous time nearest neighbor simple random walk on $B_L$ stopped upon reaching a boundary vertex or height 0.
In either case the state $a$ at the terminal vertex is known and we set $V_t(v)=a$. 

Such a system of ``state genealogy walks'' from all the vertices at time $t$ following backward in time is a dual model and is distributed as a system of {\em coalescing} simple symmetric continuous time random walks on the triangular lattice -- see for example \cite{G79}.  Since $B_L$ is finite, if $t$ is large enough all the walks starting then will with high probability hit the boundary before reaching height 0.  Indeed, if we continue the time lines and Poisson clocks below height 0 (and do not terminate the walks at height 0) then almost surely from any height $t$ there will be a random height $T_t\in (-\infty,t)$ at which the walks started from all vertices $v\in \iB$ at height $t$ will have reached boundary vertices.  What happens on the time lines below height $T_t$ does not affect $\{V_s\}_{s\ge t}$ since the states of the boundary vertices are fixed for all time.  This is equivalent to saying that the voter model itself reaches stationarity by a random finite time (distributed as $t-T_t$).  

%irrelevant to determining $V_t$ since the states of the boundary vertices are fixed.
%This defines in a unique way all the possible subsequence limiting distributions and, therefore, a unique stationary 
%distribution for $V_{t}$. 
Therefore to sample from $V_{\infty}^L$ it is enough to follow a system of coalescing continuous time simple random walks from each vertex $v$ of $\mathring{B_{L}}$ until they hit a boundary vertex $x_v$, and set $V_{\infty}^L(v)=V_{\infty}^L(x_v)$, i.e.~$V_{\infty}^L(v)=0$ if $x_v\in \partial_1\cup \partial_2$, and $V_{\infty}^L(v)=1$ otherwise.  One could instead sample from $V_{\infty}^L$ by setting $V_0(v)=2$ for every $v \in \iB$ (so the state space would become $\{0,1,2\}^{\iB}$) and simulating the voter model dynamics until there is no vertex with state 2.

\begin{figure}
\begin{center}
\includegraphics[scale=.5]{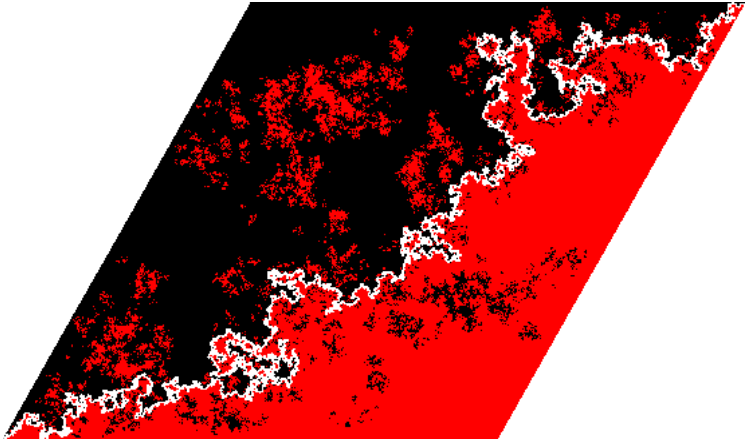}
\caption{A single realization of $V_{\infty}^L$ with $L=1026$, simulated using C++. The white curve is the exploration path or chordal interface.}
\label{fig:1sample}
\end{center}
\end{figure}

Figure \ref{fig:1sample} shows a simulation of $V_{\infty}^L$ with $L=1026$, obtained by simulating coalescing random walks from each vertex in the interior, until each one has reached the boundary.

%We have chosen the lattice and the boundary so that every configuration has a well defined interface %between clusters  of colors 0 and 1 attached to the same color boundary.  See Figure %\ref{fig:1sample}.
%By cluster we mean here a collection of vertices of the same color connected to each other by a nearest neighbor path of the same color. 

%The interface is a random curve in a dual hexagonal lattice that connects the two opposite corners of the box where the boundaries of different colors meet.  It exists for any model (dynamic or static) having the same state-space and boundary conditions.  For example, this interface curve has been studied in the context of the critical site percolation model on the triangular lattice, and has been shown to converge to SLE$_{6}$.   We can in fact describe a family of processes with the voter model and percolation as the extremes, by introducing a noise parameter into the voter model as follows.  Set $0 \leq p \leq 1$, and at each time that a Poisson clock rings, with probability $1-p$ choose a neighbor uniformly at random to adopt its color (as before), and otherwise (with probability $p$) choose a color by flipping a fair coin independent of everything else.  The original voter model corresponds to $p=0$, while for $p=1$ we recover the percolation model, as all the vertices will get their colors from iid 
%noise events. 

\subsection{Harmonic percolation and related models}
\label{sec:models}
The duality discussed in the previous section tells us that $\P(V_{\infty}(x)=1)$ (we drop the superscript $L$ when there is no ambiguity) is the probability that a simple random walk started at $x$ first hits the boundary at a 1-site.  In other words, the one-dimensional distributions of our voter model on $B_L$ are equal to those of a model we would like to call {\em harmonic percolation}. 
This is a model under which the states $\{\omega_v\}_{v \in \iB}$ are independent of each other, and as we have already suggested, $\P(\omega_v=1)$ is equal to the probability that a simple random walk started from $v$ first hits the boundary $\partial B_L$ at a 1-site (i.e., in $\partial_3\cup\partial_4$).
Harmonic percolation on an infinite strip of thickness $L$ coincides with an independent percolation model called gradient percolation \cite{Nol08,Nol09,Nol09_2}. In the case of gradient percolation the probability $p(x)$ of a site $x$ being open changes linearly from one boundary where it is 0 to the other boundary where it is 1. Thus the function $p(x)$ is harmonic inside the strip (with specified boundary conditions).
The difference between the voter and harmonic percolation models arises from the fact that the walks in the former are coalescing, whereas in the latter they are independent.  To be more explicit, coalescence in the voter model leads to non-zero correlations as in the following simple lemma.

\begin{Lemma}
\label{lem:correlation}
For any $n>1$, and any $x_1,\dots, x_n$ in the interior of $B_L$, 
\begin{align}
\P\left(\cap_{i=1}^n \{V^L_{\infty}(x_i)=1\}\right)> \prod_{i=1}^n \P\left(V^L_{\infty}(x_i)=1\right).
\end{align}
\end{Lemma}
\proof Fix $L$, $x_1$ and $x_2$, and let $\partial^0$ and $\partial^1$ denote the elements of $\partial B_L$ with fixed states 0 and 1 respectively.
%We can also see that votes of two different sites are positively correlated in the sense that for any two vertices $x$ and $y$ inside of the box %$B_{L}$,
%\begin{align}
%\P(X=1, Y=1)\ge \P(X=1)\P(Y=1)
%\end{align}
%Where $X=V_{\infty}(x)$ and $Y=V_{\infty}(y)$.
Let $S_1$ and $S_2$ be two independent random walks starting from $x_1$ and $x_2$ respectively. Let $\tau^{x_1 x_2}=\inf\{t:S_1(t)=S_2(t)\}$ be the first time that $S_1$ and $S_2$ meet each other.  Let $S'_1 = S_1$ for all times and define
\begin{align}
S'_2(t) = \begin{cases}
S_2(t), & \text{ if }t\le \tau^{x_1 x_2}\\
S_1(t), & \text{ if }t> \tau^{x_1 x_2},
\end{cases}
\end{align}
so that $S'_1$ and $S'_2$ are coalescing walks started from $x_1$ and $x_2$ respectively.  

Let $\tau^i_L=\inf\{t:S_i(t)\in \partial B_L\}$ and ${\tau'}^i_L=\inf\{t:S'_i(t)\in \partial B_L\}$ denote the respective hitting times of the boundary, and note that ${\tau'}^1_L=\tau^1_L$. Then 
\begin{align}
\P(V_{\infty}(x_1)=1,V_{\infty}(x_2)=1)=&\P(S'_1({\tau'}^1_L)\in \partial^1 ,\ S'_2({\tau'}^2_L)\in \partial^1)\nonumber\\
=&\P(S_1(\tau^1_L)\in \partial^1,\ S_2(\tau^2_L)\in \partial^1,\tau^{x_1 x_2}<\tau^1_L\wedge \tau^2_L)\nonumber\\
%\P(R_{x}\text{ hits } 1\text{ boundary},\ R_{y}\text{ hits any boundary}, R_{x}\text{ and }R_{y}\text{ meet})\\
&+\P(S_1(\tau^1_L)\in \partial^1,\ S_2(\tau^2_L)\in \partial^0,\tau^{x_1 x_2}<\tau^1_L\wedge \tau^2_L)\nonumber\\
&+\P(S_1(\tau^1_L)\in \partial^1 ,\ S_2(\tau^2_L)\in \partial^1,\tau^{x_1 x_2}\ge \tau^1_L\wedge \tau^2_L)\nonumber\\
=&\P(S_1(\tau^1_L)\in \partial^1 ,\ S_2(\tau^2_L)\in \partial^1)+\P(S_1(\tau^1_L)\in \partial^1,\ S_2(\tau^2_L)\in \partial^0,\tau^{x_1 x_2}<\tau^1_L\wedge \tau^2_L)\nonumber\\
%&+\P(R_{x}\text{ hits } 1\text{ boundary},\ R_{y}\text{ hits } 1\text{ boundary }, R_{x}\text{ and }R_{y}\text{ do not meet})\\
=&\P(V_{\infty}(x_1)=1)\P(V_{\infty}(x_2)=1)+\P(S_1(\tau^1_L)\in \partial^1,\ S_2(\tau^2_L)\in \partial^0,\tau^{x_1 x_2}<\tau^1_L\wedge \tau^2_L)\label{cor1}\\
>& \P(V_{\infty}(x_1)=1)\P(V_{\infty}(x_2)=1).\nonumber
\end{align}
This proves the result for $n=2$.  A similar coupling argument can be made for any number $n$ of walkers starting from vertices $x_{1}, \dots, x_{n} \in B_{L}$ (choosing the lower indexed random walker to continue when any two meet), establishing the claim.
%$\P(v_{1}=1,... ,v_{k}=1) \geq \P(v_{1}=1)\times ...\times \P(v_{k}=1)$. 
\Qed

For any $\epsilon>0$, if $x_1(L)$ and $x_2(L)$ are distance at least $\epsilon L$ from each other and the boundary $\del B_L$ then there exist $c_\epsilon>0$ and $C_{\epsilon}<1$ such that  $\P(V_{\infty}(x_i)=1)\in (c_{\epsilon},C_{\epsilon})$ for $i=1,2$ and all $L$, while
\[\P(S_1(\tau^1_L)\in \partial^1,\ S_2(\tau^2_L)\in \partial^0,\tau^{x_1 x_2}<\tau^1_L\wedge \tau^2_L)\le \P(\tau^{x_1 x_2}<\tau^1_L\wedge \tau^2_L)\le \P(\tau^{\Delta}_o<\tau^{\Delta}_{\partial B_{2L}}) = O\left(\frac{1}{\log{L}}\right),\]
where $\tau^{\Delta}_o$ and $\tau^{\Delta}_{\partial B_{2L}}$ are times when the difference random walk $S_1(t)-S_2(t)$ started at $x_1-x_2$ first hits the origin and the boundary of the box $B_{2L}$ respectively, and the last equality follows from Proposition 6.4.3 of \cite{LawLim10}.  Then \eqref{cor1} implies that the correlation $\rho(V_{\infty}(x_1),V_{\infty}(x_2))$ between the votes at $x_1$ and $x_2$ goes to zero as per the following.
\begin{Lemma}
\label{lem:corgoeszero}
Let $\epsilon>0$, and $x_1(L)$ and $x_2(L)$ be distance at least $\epsilon L$ from each other and the boundary $\del B_L$.  Then $\rho(V_{\infty}(x_1),V_{\infty}(x_2))\rightarrow 0$ as $L\rightarrow \infty$.
\end{Lemma} 
One can consider i.i.d.~percolation, harmonic percolation, and the stationary voter model on $B_L$ as special cases of a general 2-parameter family of models as follows.
%One could define the ordinary percolation model via random walks by starting an independent random walk at each site $v\in \iB$, and when it hits the boundary toss a fair coin (independently for each $v$) to determine the state $\omega_v$. In this way we can think of harmonic percolation as differing from percolation due to a {\em boundary coin} effect - in percolation we toss a fair coin at each point on the boundary, and in harmonic percolation the coins are completely biased. Voter model can also be coupled with the two independent percolation models. One description of a general model is the following. 
Start a continuous-time walker from each site.  Each walker initially wears a hat. Two walkers wearing hats coalesce when they meet, and instantly become a single walker wearing a hat.  Walkers not wearing hats do not coalesce with any other walkers. In addition a Poisson clock is assigned to each walker.  When such a clock rings, the walker takes a random walk step, but before doing so removes her coalescence hat with probability $q$. If a walker wearing a hat steps into a site with another walker with a hat on, the walker that just made its step becomes part of the coalescence set of the walker that was already at the site. Upon hitting a boundary site, with probability $1-p$ a walker (and its entire coalescence set) is assigned the vote of the boundary vertex it hit, and otherwise (i.e., with probability $p$) its entire coalescence set attains an independently and uniformly chosen vote.  Varying the boundary and coalescence noise parameters $p$ and $q$ between 0 and 1 allows us to interpolate between the 
four corner models:  the voter model $(p,q)=(0,0)$; harmonic percolation $(p,q)=(0,1)$; i.i.d.~percolation $(p,q)=(1,1)$; and the case $p=1,\ q=0$ corresponds to a model we would like to call cow (coalescing walk) percolation. 
%Varying these two parameters $p$--boundary noise and $q$--coalescence noise between 0 and 1 allows us to interpolate between the four extreme corner models: $p=0,\ q=0$ corresponds to the voter model, $p=1,\ q=0$ corresponds to a model we would like to call cow (coalescing walk) percolation, $p=0,\ q=1$ corresponds to the harmonic percolation model, and $p=1,\ q=1$ corresponds to the independent percolation model.
%Figure \ref{fig:pqr} shows a diagram indicating the relationships between the models.
% \begin{figure}
% \begin{center}
% \includegraphics[scale=.3]{IntroducingNoise.eps}
% \caption{{\bf may need to change this now} Introducing noise into the model allows to couple our voter %model with range of other models }
% \label{fig:pqr}
% \end{center}
% \end{figure}
 
%The cow percolation model has the same one-dimensional distributions as ordinary percolation, and the difference comes from coalescing walks replacing independent walks.  The voter model shares a similar relationship with the harmonic percolation model.

%So far this coupling has not been useful in comparing the mentioned models, other than may be suggest lemma \ref{lem:correlation}. But it seemed interesting enough to include.  

%Note that any of the parameters $p,q,r$ above could be made to depend on $L$, $v,x$ and time $t$. 

\section{Interface length}
\label{sec:length}
Recall that $|Z_L|$ denotes the length of the interface.  Since this path is a nearest neighbor simple path, there exist $c,C>0$ such that $cL\le |Z_L|\le CL^2$ almost surely.  We conjecture that 
\begin{align}
\label{lengthconj}
H_L\equiv \E[|Z_L|]\approx L^{d}
\end{align}
for some $d\in [1,2]$.  
In the case of critical i.i.d.~percolation, \eqref{lengthconj}  holds with $d=7/4=1.75$ which is also the Hausdorff dimension of the limiting law (i.e., of SLE$_6$, see \cite{Beffara08}).

For gradient percolation on an infinite strip, the interface curve between the occupied cluster and empty cluster is a.s.~unique and has expected length approximately $L^{3/7} l_{L}$, where $l_{L}$ is the horizontal length of the piece of strip in which we measure boundary length \cite[Proposition 11]{Nol08}. So, for any $\epsilon>0$, for all sufficiently large $L$, if we take a piece of strip which is $L$ long (and $L$ thick), the expected length of the interface curve $H_L$ satisfies $L^{10/7-\epsilon} \leq H_L \leq L^{10/7+\epsilon}$.
%for any positive $\epsilon$ when $L$ is large enough. 
For any $\delta > 0$, with probability going to 1 with $L$, the curve stays in the central band (around the central $L/2$ line where $p=\frac{1}{2}$) of width $L^{4/7+\delta}$ 
%where $L$ is width of the strip 
\cite[Theorem 6]{Nol08}.  Thus, as $L\rightarrow \infty$, unless we appropriately zoom in around the central line, we expect to see the rescaled interface curve converge to a straight line in the center.  Since the harmonic function inside a rhombic area with our boundary condition looks almost linear along the diagonal that connects the middle corner of the 1 valued boundary to the middle corner of the 0 valued boundary (or indeed along any parallel line), we expect that the interface curve for harmonic percolation inside our rhombus  should scale to a straight line as well.

%{\bf Check: }For gradient percolation on an infinite strip (measured inside a section of the strip of length $L$) it seems that \eqref{lengthconj} holds with $d=10/7=1.428571\dots$ (Proposition 11 \cite{Nol08}), but the limiting path is a straight line (Hausdorff dimension 1).  If one rescales properly ($\sim L^{4/7}$) one seem to get the same relation as in critical percolation case.

Writing $H_L=L^d\ell(L)$ for some function $\ell(L)$ which makes the equality true we have that 
\begin{align}
d=&\frac{\log(H_L)}{\log(L)}-\frac{\log(\ell(L))}{\log(L)}=\log_2\left(\frac{H_{2L}}{H_L}\right)-\log_2\left(\frac{\ell(2L)}{\ell(L)}\right).
\label{dformulas}
\end{align}
%and
%\begin{align}
%d=&\log_2\left(\frac{H_{2L}}{H_L}\right)+\log_2\left(\frac{\ell(L)}{\ell(2L)}\right)\label{dformulas2}.
%\end{align}
Computing the average interface curve length $\bar{Z}_m(L)$ from $m$ independent realizations of $V_{\infty}^L$ we obtain the following estimators
%As a first step we estimate the fractal dimension of the interface curve by simulations (for various noise parameters $p$). For each $L$, we sample $m=10000$ realizations of $V_{\infty}^L$ by running coalescing random walks in $B_{L}$ until they all hit the boundary, we then compute the average interface curve length $\bar{X}_{m}(L)$. 
 for $d$ based on \eqref{dformulas} 
 %and \eqref{dformulas2} respectively:
\begin{align}
\tilde{d}=\tilde{d}_{m,L}=&\frac{\log(\bar{Z}_m(L))}{\log(L)},\label{dtilde}\\
\hat{d}=\hat{d}_{m,L} =& \log_2\left(\frac{\bar{Z}_{m}(2L)}{\bar{Z}_m(L)}\right).\label{dhat}
%=\log_{2}\frac{Average\ length\ of\ curve\ in\ a\ box\ of\ size\ 2n}{Average\ length\ of\ curve\ in\ a\ box\ of\ size\ }.
\end{align}
We say that an estimator $\hat{\beta}$ (more precisely a family of estimators $\{\hat{\beta}_{m,L}:m,L \in \N\}$) is a {\em consistent} estimator of some quantity $\beta$ if 
\begin{align}
\label{consistent}
\lim_{L \ra \infty} \lim_{m \ra \infty} \hat{\beta}_{m,L}=\beta, \quad \text{almost surely}.
\end{align}
It is easy to show that $\tilde{d}$ is a consistent estimator of $d$ if and only if $\log(\ell(L))/\log(L)\rightarrow 0$ as $L\ra \infty$ (i.e., if and only if $H_L\approx L^d$), while $\hat{d}$ is a consistent estimator for $d$ if and only if $\ell(2L)/\ell(L)\rightarrow 1$ as $L\rightarrow \infty$.  Thus both estimators are consistent if $\ell$ is slowly varying at $\infty$.

If we are willing to assume that the random interface length $|Z_L|$ in a box of size $L$ satisfies $|Z_L|=CL^de^{\vep}$, where $\E[\vep]=0$, then it is natural to consider the ordinary least squares estimator $d^*$ for the slope coefficient $d$ of the simple linear regression model
\begin{align}
\log(|Z_i|)=d\log(L_i)+a+\vep_i,\label{dregression}
\end{align}
where $\{|Z_i|\}_{i\le n}$ are interface lengths on boxes of side lengths $\{L_i\}_{i \le n}$, and the $\vep_i$ are random variables with mean 0.   Note that $\ell(L)$ is constant under this assumption.
%Under this rather strong assumption

%Moreover, if $H_L=L^d\ell(L)$ where $\ell(L)$ is approximately constant for large $L$, then $\log(H_L)\approx d\log(L)+\log(C)$ and one can estimate $d$ by the slope coefficient  (i.e.~with response variable $\log(H_L)$ and explanatory variable $\log(L)$).  

The results from independent simulations of computing the average lengths of the interface curve and estimates $\hat{d},\tilde{d},d^*$ (with $m=10000$) appear in Table \ref{tab:d_estimates}.  For ordinary percolation and harmonic percolation the values are known (or expected) to be $7/4=1.75$ and $10/7\approx 1.4286$, so for these models 
%Here we mention that using the same estimators for the models where $d$ is known (ordinary percolation and gradient percolation) it seems that
 $\hat{d}$ (with $L=512$) appears to do best. For the voter model the value of $\hat{d}$ is about 1.46.

\begin{table}
\centering
\begin{tabular}{|l|l|l|l||l|l|l|l||r|}
\hline
 &   & $\hat{d}$  & & & $\tilde{d}$   &&&$d^*$\\
%$p$
$L$ & $128$ & $256$ & $512$ & $128$  &$256$ &$512$  &$1024$ &\\
\hline
% 0.0
voter &1.4427&1.4581&1.4633& 1.5947 & 1.5767& 1.5638& 1.5536 & 1.4586 \small{(s.e.~0.0012)}\\
%\hline
% 0.1&1.7493&1.7445&1.7553& & & &\\
%\hline
% 0.2&1.7446&1.7470&1.7541& & & &\\
%\hline
% 0.3&1.7444&1.7432&1.7537& & & &\\
%\hline
%0.4&1.7423&1.7448&1.7551& & & &\\
%\hline
% 0.5&1.7378&1.7455&1.7568& & & &\\
%\hline
%0.6&1.7300&1.7480&1.7515& & & &\\
%\hline
% 0.7&1.7411&1.7361&1.7575& & & &\\
%\hline
% 0.8&1.7409&1.7462&1.7517& & & &\\
%\hline
% 0.9&1.7403&1.7566&1.7421& & & &\\
%\hline
% 1.0
cow. &1.4809 &1.4876 &1.4953 &1.6689& 1.6454 &1.6279 &1.6146   &1.4878   \small{(s.e.~0.0019)}\\
harm. &1.4229 &1.4221 &1.4290& 1.6254 &1.6001 &1.5803 &1.5652  &1.4264 \small{(s.e.~0.0006)} \\
perc. &1.7401 &1.7463 & 1.7505& 1.8492 &1.8355& 1.8256 &1.8181 & 1.7458 \small{(s.e.~0.0019)}\\
\hline
\end{tabular}
\caption{Estimates of $d$ with $m=10000$ and 
%for various values of $p$, 
%computed from \eqref{dhat} 
with $L=128,256,512$ (and $1024$), rounded to 4 decimal places. For harmonic percolation and percolation, $\hat{d}$ is the closest estimator to the true values of 1.42857 and 1.75 respectively.  }
\label{tab:d_estimates}
\end{table}
%Proceeding in this way, the estimated value of $1.7535$ for $d$ (percolation model) based on %\eqref{dhat} with $L=512$ is fairly close to the true limiting value of 1.75. The corresponding %estimate for the voter model is around 1.45 (see Table \ref{tab:d_estimates}). 
%An interesting observation is that as soon as we introduce positive but fixed noise, the value for the dimension jumps to around 1.75 which is the value for percolation model.  {\bf Can we prove anything like this by e.g.~renormalization?}

\section{The sizes of coalescing classes and related questions}
\label{sec:coalescing}
The interface curve cannot pass through any connected cluster of common votes.  
%In the voter and cow percolation models, two vertices certainly share the same vote/state if the walks started from them coalesce before hitting the boundary. 
The difference between the voter and harmonic percolation models is that the states are determined by coalescing random walks rather than independent random walks (started at each site).  If the coalescing classes in the voter model are negligible as $L \uparrow \infty$, both in terms of size and the correlation between votes in different classes, then perhaps some kind of rescaling argument would allow one to compare the voter model to the harmonic percolation model.  One expects that the rescaled interface curve for harmonic percolation on $B_L$ converges to a straight line (Pierre Nolin has proved this on the strip \cite{Nol08}), so one might expect the same to be true for the voter model, if the coalescing classes are indeed negligible as $L \uparrow \infty$. 

Clustering behaviour for the 2-dimensional voter model has been well studied in the probability literature (see e.g.~\cite{CG86}), but (as far as we know) not in the current setting of a finite domain with unflinching boundary. 
As a small step in the direction of understanding the correlation between votes in different classes, let us verify that any two sites $x$ and $y$ are less likely to share a common vote (than they otherwise would be) if they are not in the same coalescing class.  Fix $L$ and start coalescing walks from every site in $\mathring{B_{L}}$.  The walks define an equivalence relation on $\mathring{B_{L}}$ in the sense that $x\sim y$ if and only if the walks started from $x$ and $y$ coalesce before hitting the boundary.  Let $C_L(x)$ denote the (random) equivalence class of $x$.
%(i.e.~those vertices such that the walks started from them coalesced with the walk from $x$ before hitting the boundary). 
Let $A=A_{xy}=\{V_{\infty}(x)=V_{\infty}(y)\}$ and $B=B_{xy}=\{C_L(x)=C_L(y)\}=\{x\sim y\}$.  Then $\P(B)\in (0,1)$ and
\begin{align*}
\P(A)=\P(A|B)\P(B)+\P(A|B^c)\P(B^c)=\P(B)+\P(A|B^c)\P(B^c)>\P(A|B^c),
\end{align*}
as claimed.  On the other hand, as we have seen earlier, if $x=x(L)$ and $y=y(L)$ are distance at least $\epsilon L$ apart then $\P(B^c)\rightarrow 1$ and $\P(A|B^c)-\P(A)\rightarrow 0$ as $L\nearrow \infty$, with $\P(A)$ being bounded away from zero as $L\nearrow \infty$ if $x$ and $y$ are also at least $\epsilon L$ distance from the boundary.

\begin{figure}
\begin{center}
\includegraphics[width=346pt,height=200pt]{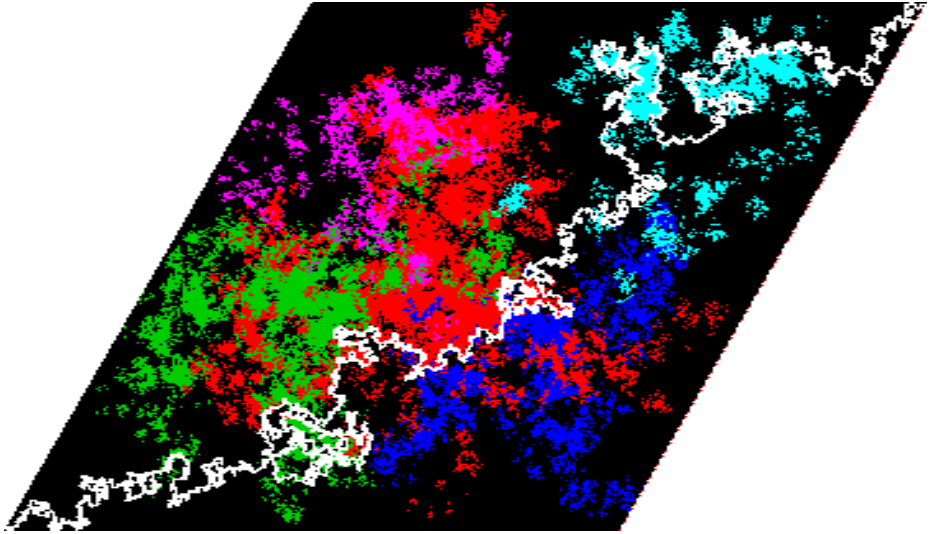}
\caption{A single realization of the 5 largest coalescing classes (colors other than black) for $L=1026$. The white curve is the chordal interface.}
\label{fig:5largest}
\end{center}
\end{figure}

We are hereafter interested in 
the behavior of the expected size of the class of the centre of the rhombus $B_L$, which we will for convenience take to be the origin (if $L$ is not odd we consider the centre/origin to be any one of the closest vertices to the centre)
%(i.e., we consider the shifted rhombic box of side $L$ centred at the origin)
and the expected size of the largest  class $\E\left[|M_L|\right]$ as $L\ra \infty$, where
\begin{align}
\label{ML}
 M_L=C_L(x'), &\text{ where $x'$ is chosen such that } |C_L(x')|=\max_{x \in \mathring{B_{L}}}|C_L(x)|,
\end{align} 
where some tie-breaking rule is used to choose $x$, if necessary.  In particular, we ask what proportion of all vertices in the box are in the largest class, as $L\ra \infty$? Since $|\mathring{B_{L}}|\approx c L^2$, we are interested in $\lim_{L\ra \infty}\E
\left[|M_L|/L^2\right]$.  Figure \ref{fig:5largest} shows a single realization of the 5 largest classes for $L=1026$.

The following is our main rigorous result, which shows that coalescing classses have (on average) small (but only logarithmically small) volume compared to the whole box.  We note that the lower bound can be improved slightly with a little more effort.
\begin{Theorem}
\label{thm:ccp_lim}
There are constants $C'$ and $C''$ in $(0,\infty)$, such that  
\begin{align}
\frac{C'}{\left(\log{L}\right)^{713}} \leq \E\left[\frac{|C_{L}(o)|}{|\mathring{B_{L}}|}\right] \leq \E\left[\frac{|M_{L}|}{|\mathring{B_{L}}|}\right]\leq \left(\frac{C''}{\log{L}}\right)^{1/2}.\label{thmbounds}
\end{align}
\end{Theorem}
\proof
First note that to have $B_L$ centered at the origin and to have the smaller rhombi used in the proof to be consistent with the lattice we need $L$ to be such that $L-1$ is divisible by 12, and then for rhombi of fractional side lengths such as $L/2$ we should use $B_{\frac{L-1}{2}+1}$ instead of $B_{L/2}$.  For notational convenience we will ignore these issues, but we note though that the same arguments would in any case work with trivial but messy modifications. 

For both the upper and lower bounds in \eqref{thmbounds} we will use the fact that, for $x \in \mathring{B_{L}}$,
\begin{align}
\label{equality}
 \E\left[|C_{L}(x)|\right] &= \E\left[\sum_{y \in \mathring{B_{L}}} \indic{y \in C_{L}(x)}\right]= \sum_{y \in \mathring{B_{L}}}\P(y \in C_{L}(x)).
 \end{align}

To verify the lower bound, let $\{S_x\}_{x\in \mathring{B_{L}}}=\{\{S_x(t)\}_{t\ge 0}\}_{x\in \mathring{B_{L}}}$ be independent continuous-time (with jump rate 1), nearest-neighbor random walks on the triangular lattice, with respective starting points $S_x(0)=x\in \mathring{B_{L}}$.  For $x,y \in \mathring{B_{L}}$, let $\tau^{xy}=\inf\{t:S_x(t)=S_y(t)\}$, and $\tau^x_{L}=\inf\{t:S_x(t)\in \partial B_L\}$ denote the meeting times and boundary hitting times respectively. 
Then \eqref{equality} with $x=o$ can be written as
\begin{align}
 \label{lowBoundCoalClass1}
  \E\left[|C_L(o)|\right] = \sum_{x \in \mathring{B_{L}}} \P\left( \tau^{ox} < \tau^o_{L} \wedge \tau^x_{L} \right). 
\end{align}
By Lemma \ref{lemma_a} below there is positive constant $c'$ such that for each $x \in B_{L/2} \setminus B_{L/4}$,  
\begin{align}
 \P\left( \tau^{ox} < \tau^o_{L} \wedge \tau^x_{L} \right) \geq \frac{c'}{(\log L)^{713}}\label{lem722}.
\end{align}
Combining \eqref{lowBoundCoalClass1} and \eqref{lem722} we obtain
\begin{align}
 \E\left[\frac{|C_L(o)|}{|\mathring{B_{L}}|}\right] &\geq \frac{1}{cL^2}\sum_{x \in B_{L/2} \setminus B_{L/4}} \P\left( \tau^{ox} < \tau^o_{L} \wedge \tau^x_{L} \right)\geq \frac{C'}{(\log L)^{713}},\label{lowerdone}
\end{align}
for another positive constant $C'$, which verifies the lower bound in \eqref{thmbounds}.

To establish the upper bound, note that for any $x,y\in \mathring{B}_L$ the difference walk $S^{\Delta}_{x-y}(t) = S_x(t)-S_y(t)$ is also a simple symmetric random walk started at $x-y$ but with jump rate 2. Let $\tau^{\Delta_x}_o = \inf\{t:S^{\Delta}_x(t)=o\}$ and $\tau^{\Delta_x}_{L}= \inf\{t:S^{\Delta}_x(t) \in \partial B_L\}$ be the first hitting times of the origin and the boundary $\partial B_{L}$ by the difference walk $S^{\Delta}_x$, and let $\tau^{x}_o = \inf\{t:S_x(t)=o\}$ be the first time $S_x$ hits the origin. Then 
\begin{align}
 \frac{\E[|C_L(x)|]}{cL^2} &= \frac{1}{cL^2}\sum_{y \in \mathring{B}_L} \P(\tau^{xy}<\tau^{x}_{L} \wedge \tau^{y}_{L})\nonumber\\
 &\leq \frac{1}{cL^2}\sum_{y \in \mathring{B}_{L}} \P(\tau^{\Delta_{x-y}}_{o}<\tau^{\Delta_{x-y}}_{2L})\nonumber\\
 &\leq \frac{1}{cL^2}\sum_{y \in B_{2L}\setminus o} \P(\tau^{y}_{o} < \tau^{y}_{2L}) + \frac{\P(\tau^{o}_{o} < \tau^{o}_{2L})}{cL^2}\label{bananaman1}.
 \end{align}
Using Theorem 6.4.3 of \cite{LawLim10} on the summation, \eqref{bananaman1} is bounded above by
\begin{align}
\frac{C}{L^2}\sum_{y \in B_{2L}\setminus o}\frac{\log(L/|y|)}{\log(L)}+\frac{1}{cL^2}.
% - \log (c'|y|) +o(|y|^{-1}) + O\left( \frac{1}{\log(L)}\right)}{\log(L)}\right]+\frac{1}{cL^2}\nonumber 
% & \frac{1}{cL^2}\sum_{y \in B_{2L}\setminus o}\left[ \frac{\log(L) - \log (c'|y|) +o(|y|^{-1}) + O\left( \frac{1}{\log(L)}\right)}{\log(L)}\right]+\frac{1}{cL^2}\nonumber \\
% &\leq \frac{1}{cL^2}\sum_{y \in B_{2L}\setminus o} \frac{\log(L) - \log (c'|y|) + 1}{\log(L)} + O\left( \left(\frac{1}{\log(L)}\right)^2\right).
\label{limlaw}
 \end{align}
Next, we split the sum into dyadic annuli (all but finitely many of which contain no vertices). Since $B_L$ has been defined via the number of vertices on the boundary (so only for integer $L$), we let $\hat{B}_R$ denote the rhombic $R\times R$ box centered at the origin in $\R^2$ and let $A_{L,k}$ denote the (possibly empty) intersection of ($\hat{B}_{2L/2^{k-1}}\setminus \hat{B}_{2L/2^{k}}$) with the triangular lattice. Then 
the first term in \eqref{limlaw} is equal to
 \begin{align}
\frac{C}{L^2}\sum_{k=1}^{\infty}\ \sum_{y \in A_{L,k}}\frac{\log(L/|y|)}{\log(L)}\le & \frac{C}{L^2}\sum_{k=1}^{\infty}\ \sum_{y \in A_{L,k}}\frac{\log(2^{k-1})}{\log(L)}\nonumber\\
\le & \frac{C'}{L^2}\sum_{k=1}^{\infty}\left(\frac{2L}{2^k}\right)^2\cdot \frac{k-1}{\log L}\nonumber\\
% & \frac{1}{cL^2}\sum_{k=1}^{\infty}\ \sum_{y \in B_{2L/2^{k-1}}\setminus B_{2L/2^{k}}}\frac{\log(L) - \log (c'|y|) + 1}{\log(L)} + O\left( \left(\frac{1}{\log(L)}\right)^2\right)\nonumber\\
% &\leq \frac{1}{cL^2}\sum_{k=1}^{\infty}\ \sum_{y \in B_{2L/2^{k-1}}\setminus B_{2L/2^{k}}}\frac{\log(L) - \log \left(c'\frac{L}{2^{k}}\right) + 1}{\log(L)} + O\left( \left(\frac{1}{\log(L)}\right)^2\right)\nonumber\\
%&\leq\frac{1}{cL^2}\sum_{k=1}^{\infty}\ \tilde{C}\left(\frac{2L}{2^k}\right)^2\frac{\log (L) - \log (L) + k\log 2 + c''}{\log(L)} + O\left( \left(\frac{1}{\log(L)}\right)^2\right)\nonumber\\
%=\frac{C'}{\log(L)}\sum_{k=1}^{\infty} \frac{k\log2 + c''}{4^k} +  O\left( \left(\frac{1}{\log(L)}\right)^2\right)
\le &\frac{C''}{\log L}.\nonumber
\end{align}
Therefore
\begin{align}
 \sup_{x \in \mathring{B_{L}}}\E\left[\frac{|C_{L}(x)|}{|\mathring{B_{L}}|}\right] \leq \frac{C''' }{\log{L}}.
\end{align}
Markov's inequality gives 
$
(\epsilon |\mathring{B_{L}}|)^{-1} \E \left[|C_{L}(x)|\right] \geq \P(|C_{L}(x)|>\epsilon |\mathring{B_{L}}|)
$ for any $\epsilon=\epsilon(L)>0$,
so 
\begin{align}
\sup_{x \in \mathring{B_{L}}}\P(|C_{L}(x)|>\epsilon |\mathring{B_{L}}|) \leq \frac{1}{\epsilon}\sup_{x \in \mathring{B_{L}}}\E\left[\frac{|C_{L}(x)|}{|\mathring{B_{L}}|}\right] \leq \frac{1}{\epsilon}\cdot \frac{C''' }{\log{L}}.
\end{align}
It follows that
\begin{align}
 \E\left[\frac{|M_{L}|}{|\mathring{B_{L}}|}\right] &\leq \frac{1}{|\mathring{B_{L}}|}\E\left[|M_{L}|\indic{|M_{L}| > \epsilon|\mathring{B_{L}}|}\right] + \epsilon = \frac{1}{|\mathring{B_{L}}|}\E\left[\sum_{x\in \mathring{B_{L}}}\indic{x\in M_{L}}\indic{|M_L|>\epsilon|\mathring{B_{L}}|} \right] + \epsilon\nonumber\\
 &\leq \frac{1}{|\mathring{B_{L}}|}\E\left[\sum_{x\in \mathring{B_{L}}}\indic{|C_{L}(x)|>\epsilon|\mathring{B_{L}}|} \right] + \epsilon=\frac{1}{|\mathring{B_{L}}|}\sum_{x \in \mathring{B_{L}}} \P(|C_{L}(x)|>\epsilon|\mathring{B_{L}}|)  + \epsilon \\
 &\leq \frac{1}{|\mathring{B_{L}}|}\sum_{x\in \mathring{B_{L}}}\sup_{x\in \mathring{B_{L}}}\P(|C_{L}(x)|>\epsilon|\mathring{B_{L}}|)  + \epsilon= \sup_{x\in \mathring{B_{L}}}\P(|C_{L}(x)|>\epsilon|\mathring{B_{L}}|)  + \epsilon\nonumber\\
 & \leq \frac{1}{\epsilon}\cdot\frac{C''}{\log{L}} + \epsilon.\nonumber
\end{align}
Choose  $\epsilon(L) = (\frac{1}{\log{L}})^{1/2}$ to get the claimed upper bound. 
\Qed

\begin{Lemma}
\label{lemma_a}
Fix some small $\epsilon>0$. There exists $c'>0$, such that for all $x \in B_{L/2}\setminus B_{\epsilon L}$, 
\begin{align}
\label{Lemma3.2}
 \P\left( \tau^{ox} < \tau^o_{L} \wedge \tau^x_{L} \right) \geq \frac{c'}{(\log L)^{713}}.
\end{align}
\end{Lemma}
\proof
%Proof of Lemma a
The triangular lattice is constructed from 3 families of parallel lines (or ``directions''), denoted by $D=\{\neswarrow,\nwsearrow,\leftrightarrow\}$, with each vertex being at the intersection of 3 such lines (one from each family), and having 6 nearest neighbors corresponding to moving ``up'' or ``down'' in any one of these directions.  
We define a system of two dependent {\em discrete-time} random walks $\hat{S}_o$ and $\hat{S}_x$ on the triangular lattice in the following way: $\hat{S}_o(0)=o$ and $\hat{S}_x(0)=x$; for each $t \in \N$ toss a fair coin to decide which of $\hat{S}_o(t)$ or $\hat{S}_x(t)$ makes an i.i.d.~uniformly chosen nearest-neighbor step on the triangular lattice (while the other does not move). 
Let $\hat{\tau}^{ox}$, $\hat{\tau}^x_{L}$, and $\hat{\tau}^o_{L}$ be the corresponding meeting and boundary hitting times for $\hat{S}_o(t)$ and $\hat{S}_x(t)$.  Since $\{(\hat{S}_o(t),\hat{S}_x(t))\}_{t \in \Z_+}$ has the same law as the jump process of $\{(S_o(t),S_x(t))\}_{t\in \R_+}$ and the event in \eqref{Lemma3.2} depends only on the relative sizes of the hitting times, we have  
\begin{align}
 \P\left( \tau^{ox} < \tau^o_{L} \wedge \tau^x_{L} \right) = \P\left( \hat{\tau}^{ox} < \hat{\tau}^o_{L} \wedge \hat{\tau}^x_{L} \right).\label{pants1}
\end{align}
Now we focus on the discrete time random walks $\hat{S}_o(t)$ and $\hat{S}_x(t)$.
For $t \in \N$, let $\mathcal{R}_{t}$ denote the set of ordered partitions $r = (r_\neswarrow,r_\nwsearrow,r_\leftrightarrow)(t)$ of $\{1,2,\dots,t\}$ into three (possibly empty) sets.  For $s\le t$, let $h(s)\in D$ denote the direction of the step taken by (one of) the pair $(\hat{S}_o,\hat{S}_x)$ at time $s$.  For $r=(r_\neswarrow,r_\nwsearrow,r_\leftrightarrow)\in \mathcal{R}_t$ let $A_{r}$ be the event that for each $\bullet \in \{\neswarrow,\nwsearrow,\leftrightarrow\}$ and $s\in r_\bullet$, $h(s)=\bullet$, i.e., that steps in direction $\neswarrow$ are taken at times in $r_{\neswarrow}$ etc.
Conditioning on $A_{r}$ we can rewrite the probability above as follows
\begin{align}
\label{lowBoundCoalClass2}
 \P\left( \hat{\tau}^{ox} < \hat{\tau}^o_{L} \wedge \hat{\tau}^x_{L} \right) = \sum_{t=0}^{\infty} \sum_{r \in \mathcal{R}_{t}}\P\left(\hat{\tau}^{ox} = t, \hat{\tau}^o_{L} \wedge \hat{\tau}^x_{L}>t\ |\ A_{r} \right)\P\left(A_{r}\right).
\end{align}
Let $S^{\Delta}$ and $S^{\Sigma}$ denote the difference and sum walks starting at $x$, defined by $S^{\Delta}(t)=\hat{S}_x(t)-\hat{S}_o(t)$ and $S^{\Sigma}(t)=\hat{S}_x(t)+\hat{S}_o(t)$, and let $\tau^{\Delta}_o$ and $\tau^{\Delta}_{L}$, and $\tau^{\Sigma}_o$ and $\tau^{\Sigma}_{L}$ be the hitting times of the origin and the boundary $\partial B_{L}$ by the difference and sum walks respectively.
Since
\begin{align}
 \hat{S}_x(t)= \frac{\hat{S}_x(t)-\hat{S}_o(t)}{2}+\frac{\hat{S}_x(t)+\hat{S}_o(t)}{2},
\end{align}
we have 
\begin{align}
 |\hat{S}_x(t)|\geq L => |\hat{S}_x(t)-\hat{S}_o(t)| \geq L \text{ or } |\hat{S}_x(t)+\hat{S}_o(t)|\geq L.
\end{align}
A similar statement holds for $\hat{S}_o(t)$. Thus we have
\begin{align}
 \{\hat{\tau}^o_{L} \wedge \hat{\tau}^x_{L}>t\} \supseteq \{\tau^{\Sigma}_{L}>t\} \cap \{\tau^{\Delta}_{L} >t \}.
\end{align}
Therefore \eqref{lowBoundCoalClass2} can be continued as follows
\begin{align}
 &\sum_{t=0}^{\infty} \sum_{r\in \mathcal{R}_{t}}\P\left(\hat{\tau}^{ox} = t, \hat{\tau}^o_{L} \wedge \hat{\tau}^x_{L}>t\ |\ A_r \right)\P\left(A_{r}\right)\nonumber\\
 &\hspace{0.5cm} \geq \sum_{t=0}^{\infty} \sum_{r\in \mathcal{R}_{t}}\P\left(\tau^{\Delta}_o = t, \tau^{\Delta}_{L}>t, \tau^{\Sigma}_{L}>t\ |\ A_{r} \right)\P\left(A_{r}\right)\\
 &\hspace{0.5cm}= \sum_{t=0}^{\infty} \sum_{r\in \mathcal{R}_{t}}\P\left(\tau^{\Delta}_o = t, \tau^{\Delta}_{L}>t|A_{r}\right)\P\left( \tau^{\Sigma}_{L}>t\ |\ A_{r} \right)\P\left(A_{r}\right),\nonumber
\end{align}
where the last equality follows from the fact that the sum and difference walks are conditionally independent given $A_r$ (e.g., if we know that the sum walk makes a positive step in a specific direction, the difference walk is still equally likely to make either a positive or negative step in that direction).
%we now explain the last equality.
%As we know the direction along which each step of either $\hat{S}_x$ or $\hat{S}_o$ is taken, the sum and %difference random walks are reduced to being sums of {\em independent} discrete time one-dimensional random %walks in each of the three directions $\{\neswarrow,\nwsearrow,\leftrightarrow\}$ 

Let $x\in B_{L/2}$.   Truncating the infinite sum and using Lemma \ref{lem:projections} below we have
\begin{align}
\label{lowBoundCoalClass3}
 &\sum_{t=0}^{\infty} \sum_{r\in \mathcal{R}_{t}}\P\left(\tau^{\Delta}_o = t, \tau^{\Delta}_{L}>t|A_{r}\right)\P\left( \tau^{\Sigma}_{L}>t\ |\ A_{r} \right)\P\left(A_{r}\right)\nonumber\\
 &\hspace{0.5cm}\geq \sum_{t=0}^{n} \sum_{r\in \mathcal{R}_{t}}\P\left(\tau^{\Delta}_o = t, \tau^{\Delta}_{L}>t|A_{r}\right)\P\left( \tau^{\Sigma}_{L}>t\ |\ A_{r} \right)\P\left(A_{r}\right)\\
 &\hspace{0.5cm}\geq \P\left( \check{\tau}^{o}_{\pm L/12}>n \right)\sum_{t=0}^{n} \sum_{r\in \mathcal{R}_{t}}\P\left(\tau^{\Delta}_o = t, \tau^{\Delta}_{L}>t\ |\ A_{r}\right)\P\left(A_{r}\right),\nonumber
\end{align}
where $\check{\tau}^{u}_{\pm K}$ denotes the (discrete) time a one-dimensional simple symmetric random walk started at $u\in[-K,K]$ first hits $+K$ or $-K$.  Summarizing from \eqref{lowBoundCoalClass2} until this point, and continuing we have
 \begin{align}
 \P\left( \hat{\tau}^{ox} < \hat{\tau}^o_{L} \wedge \hat{\tau}^x_{L} \right)&\hspace{0.5cm} \ge\P\left( \check{\tau}^{o}_{\pm L/12}>n \right)\sum_{t=0}^{n} \sum_{r\in \mathcal{R}_{t}}\P\left(\tau^{\Delta}_o = t, \tau^{\Delta}_{L}>t\ |\ A_{r}\right)\P\left(A_{r}\right)\nonumber\\
 &\hspace{0.5cm} = \P\left( \check{\tau}^{o}_{\pm L/12}>n \right)\P\left( \tau^{\Delta}_o \leq n, \tau^{\Delta}_o < \tau^{\Delta}_{L}\right) \nonumber \\
 &\hspace{0.5cm} \geq \P\left( \check{\tau}^{o}_{\pm L/12}>n \right)\P\left( \tau^{\Delta}_o < \tau^{\Delta}_{L}\leq n\right)\nonumber \\
 &\hspace{0.5cm} \geq \P\left( \check{\tau}^{o}_{\pm L/12}>n \right)\left[\P\left( \tau^{\Delta}_o < \tau^{\Delta}_{L}\right) - \P\left( \tau^{\Delta}_{L} > n \right)\right]\nonumber\\
 & \hspace{0.5cm}\geq \P\left( \check{\tau}^{o}_{\pm L/12}>n \right)\left[\P\left( \tau^{\Delta}_o < \tau^{\Delta}_{L}\right) - \P\left( \check{\tau}^{o}_{\pm L/2}>\frac{n}{2} \right)\right], \label{projections}
\end{align}
%where $\tau^{1,o}_{\pm L}$ is the first time a one-dimensional simple random walk starting at the origin hits $\{L,-L\}$. 
where the last inequality follows from Lemma \ref{lem:projections} below.

Let $R_{n}=|\{z \in \Z:S^o(t) =z \text{ for some } t\leq n\}|$ be the size of the range of a one-dimensional discrete-time random walk $S^o$ (started at $o$) up to time $n$.
%:
%\begin{align}
% R_n = . 
%\end{align}
Then 
\begin{align}
\P(R_n\leq L)\le  \P(\check{\tau}^{o}_{\pm L} >n) \leq \P(R_{n} \leq 2L).\label{rangetime}
\end{align}
%and
%\begin{align}
% \P(\tau^{1,o}_{\pm L/12} >n) \geq \P(R_n\leq L/12).
%\end{align}
According to Theorem 2 of \cite{Chen06}, for any sequence $b_n= o(n)$ diverging to $+\infty$ 
\begin{align}
\label{RangeTheorem}
 \lim_{n\rightarrow \infty} \frac{1}{b_n}\log \P\left(R_n \leq \sqrt{\frac{n}{b_n}}\right) = -\frac{\pi^2\sigma}{2},
\end{align}
where $\sigma = \E\left[S^o(1)^2\right] =1$. 
Letting $n = L^2\log\log L$, and $b_n = (12)^2\log \log L$, \eqref{rangetime} and \eqref{RangeTheorem} yield (for $\delta \in (0,1/711)$) that for large $n$ (and $L$),
%for a small $0<\delta < 1/720$ and $L$ large enough we obtain
\begin{align}
\P\left( \check{\tau}^{o}_{\pm L/12}>n \right) \ge \P(R_n\leq \sqrt{n/b_n}) \geq \frac{1}{(\log L)^{\frac{(12\pi)^2(1+\delta)}{2}}} \geq \frac{1}{(\log L)^{712}}.\label{712}
\end{align}
Similarly, with $n/2=\frac{L^2}{2}\log \log L$ and $b_{n/2}= \frac{1}{2}\log \log L$, \eqref{rangetime} and \eqref{RangeTheorem} yield (for $\delta\in (0,1-8/\pi^2)$)  that for large $n$ (and $L$),
\begin{align}
\P\left( \check{\tau}^{o}_{\pm L/2}>\frac{n}{2} \right)\le \P\left(R_{n/2}\leq \sqrt{(n/2)/(b_{n/2})}\right) \leq \frac{1}{(\log L)^{\frac{\pi^2}{4}(1-\delta)}}\leq \frac{1}{(\log L)^2}.\label{1312}
\end{align}

According to Theorem 6.4.3 of \cite{LawLim10} the term $\P\left( \tau^{\Delta}_o < \tau^{\Delta}_{L}\right)$ for $x \in B_{L/2} \setminus B_{\epsilon L}$ can be bounded uniformly from below and above by $\frac{c}{\log L}$ and $\frac{C}{\log L}$ for some positive constants $c$ and $C$ and $L$ large enough.
% On the other hand, let $S^o_t$ and $B^o_t$ in the lines below indicate positions of simple symmetric random walk and Brownian motion started from the origin at time $t$, then 
% \begin{align}
% \label{hitBounndaryProb1}
%  \P\left( \tau^{1}_{aL} > bn \right) &= 1- \P\left( \tau^{1}_{aL} \leq bn \right) \nonumber \\
%  &= 1- \P\left( \sup_{0<r\leq bn} S^o_r \geq aL\right) \nonumber \\
%  &= 1-2\P\left( S^o_{bn} > aL\right) - \P\left( S^o_{bn} = aL\right)  \\
%  &= \P\left(-aL \leq S^o_{bn} < aL\right),\nonumber \\
%  &= \P\left( |B^o_1| \leq \frac{aL}{\sqrt{bn}}\right) \pm O\left( \frac{c'}{\sqrt{bn}}\cdot \frac{2aL}{\sqrt{bn}}\right),\nonumber \\
% \end{align}
% for some positive constant $c'$, with the third equality due to reflection principle and the second inequality due to Berry–Esseen theorem. 
% And therefore
% \begin{align}
%  \frac{caL}{\sqrt{bn}} \leq \P\left( \tau^{1}_{aL} > bn \right)\leq \frac{CaL}{\sqrt{bn}}, \nonumber
% \end{align}
% for some positive constants $c$ and $C$ and $L$ large enough.
Inserting this estimate, \eqref{712} and  \eqref{1312} into \eqref{projections} verifies that there exists a constant $c'>0$ such that for all $L$, uniformly in $x \in B_{L/2} \setminus B_{\epsilon L}$,
\begin{align}
\P\left( \hat{\tau}^{ox} < \hat{\tau}^o_{L} \wedge \hat{\tau}^x_{L} \right)\ge \frac{c'}{(\log L)^{713}},
\end{align}
as required.
\qed
%end of proof of lemma a

\begin{Lemma}
\label{lem:projections}
With the definitions of $\check{\tau}^{o}_{\pm L}$, $\tau^{\Sigma}_{L}$, $\tau^{\Delta}_{L}$, and $A_{r}$ as in the proof of Lemma \ref{lemma_a}, for any $r \in \mathcal{R}_{t}$ and $t\leq n$ we have 
%{\bf Give the definitions of the terms here again, and any required qualifiers on $x$ etc.}
\begin{align}
\P\left( \check{\tau}^{o}_{\pm L/12}>n \right)\le &\P\left( \tau^{\Sigma}_{L}>t\ |\ A_{r} \right)\\
\P\left( \tau^{\Delta}_{L} > n \right)  \le &\P\left( \check{\tau}^{o}_{\pm L/2}>\frac{n}{2} \right).
\end{align}
\proof 
To verify the first claim, first recall the definition of $\mathcal{R}_t$ after \eqref{pants1}.  For each $r \in \mathcal{R}_{t}$ we construct a two-dimensional random walk on the triangular lattice (started at $x \in B_{L/2}$, with $i=|r_{\neswarrow}|,j=|r_{\nwsearrow}|$ and $k=|r_{\leftrightarrow}|$ steps along the three directions in $D$ respectively) from a one-dimensional random walk of $t$ steps in the following way: let the one-dimensional walk be $S^1(t)=\sum_{\ell=1}^t X_\ell$ with each $X_\ell\in \{\pm 1\}$; designate the first $i$ steps to be in direction ``$\neswarrow$'', the next $j$ steps to be in direction ``$\nwsearrow$'', and the final $k$ steps to be in direction ``$\leftrightarrow$''; construct the two-dimensional random walk starting at $x$ by picking steps from each group according to the partition $r$ (preserving the order of steps within each of the groups). If the one-dimensional walk started at the origin stays confined to the interval $(-L/12,L/12)$, then the first $i$ steps, next 
$j$ steps and 
next $k$ steps have 
displacements from their respective starting points at most $L/12$, $2L/12$, and $2L/12$ respectively and the two-dimensional walk started at $x$ stays confined to $\mathring{B_{L}}$.  Therefore we have for $t \leq n$ that
\begin{align}
  \P\left( \tau^{\Sigma}_{L}>t\ |\ A_{r} \right) \geq \P\left( \check{\tau}^{o}_{\pm L/12}>t \right) \geq \P\left( \check{\tau}^{o}_{\pm L/12}>n \right).
\end{align}
This verifies the first claim.

\begin{figure}
\begin{center}
\includegraphics[scale=.15]{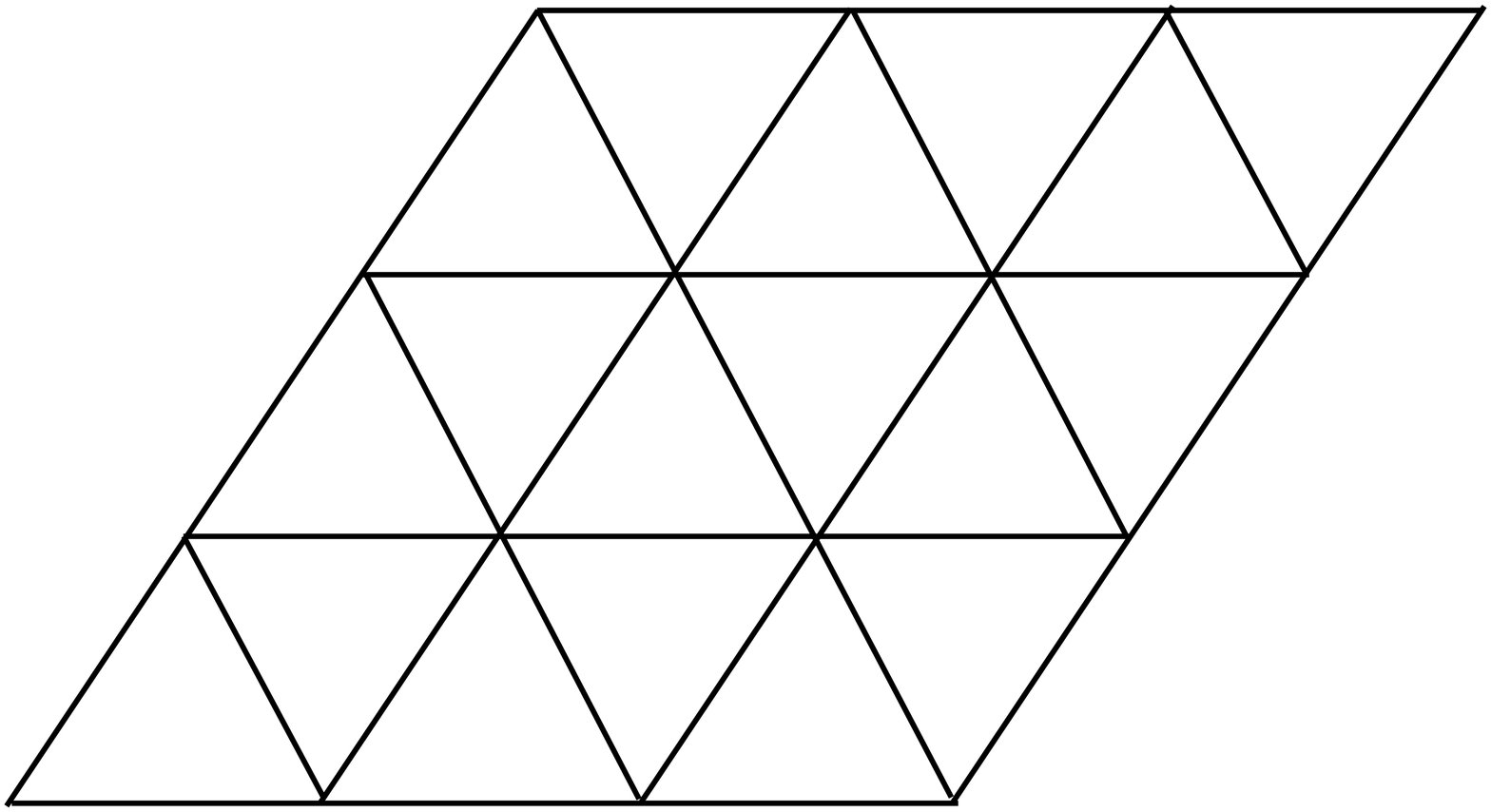}
{\Huge{$\qquad \overset{\Phi}{\longrightarrow} \qquad$}}
\includegraphics[scale=.15]{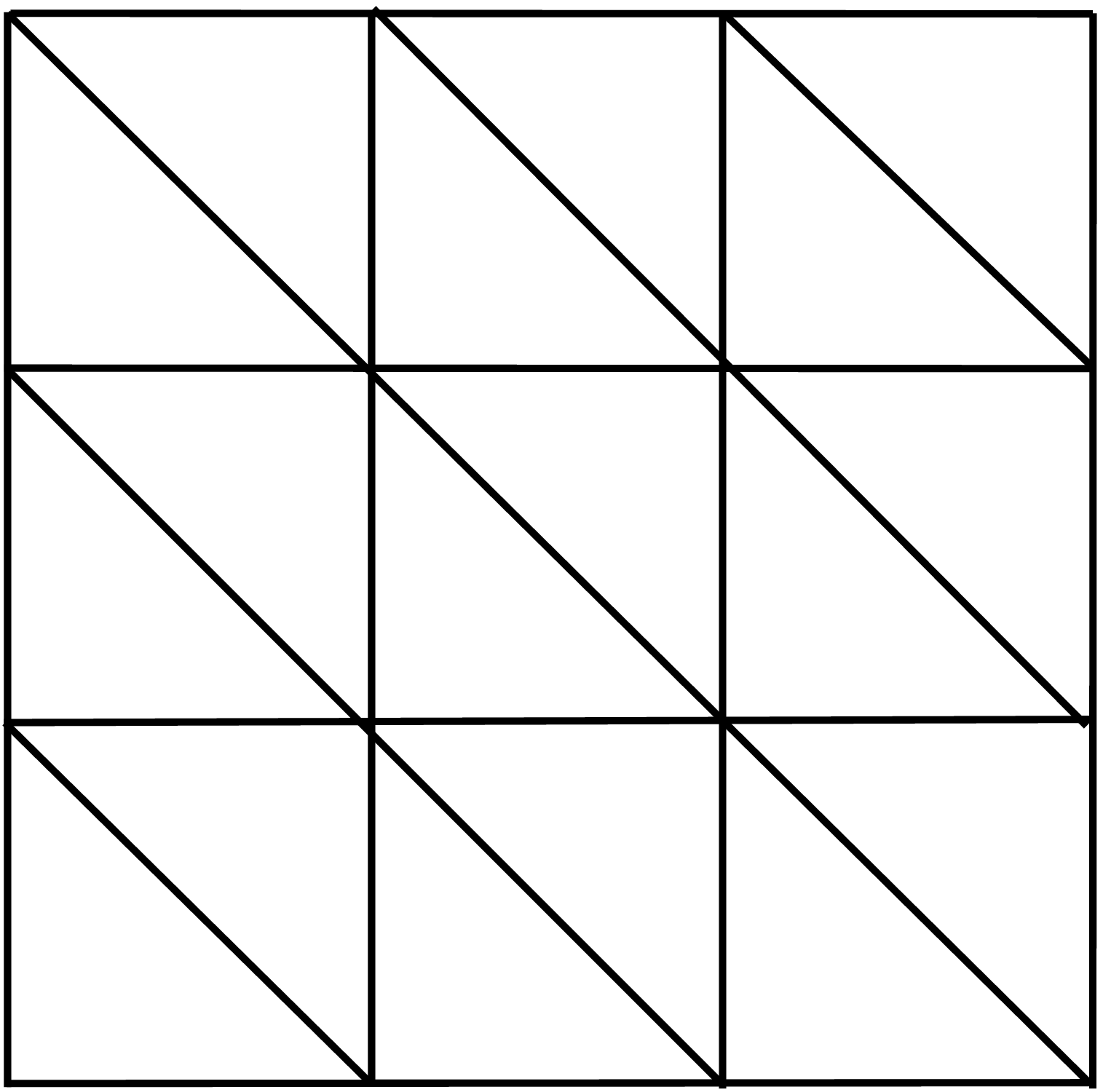}
\caption{A mapping between two planar embeddings of the triangular lattice.}
\label{fig:LinearMap}
\end{center}
\end{figure}

For the second claim, we consider two one-dimensional random walks, $S^{1,\Delta}$ and $S^{2,\Delta}$, that are the following ``projections'' of the two-dimensional discrete-time (difference) random walk $S^{\Delta}$ starting at $x$ onto the lines parallel to the two sides of the rhombic box $B_L$. Under the linear transformation $\Phi$ depicted in Figure \ref{fig:LinearMap} ``projections'' are simply standard orthogonal projections onto the two coordinate axes.
Thus, the walk $S^{1,\Delta}$ makes no step when $S^{\Delta}$ steps in the direction $\neswarrow$, makes a step -1 when the increment of $S^{\Delta}$ is either $\leftarrow$ or $\nwarrow$, and +1 when the increment of $S^{\Delta}$ is either $\rightarrow$ or $\searrow$.  Similarly $S^{2,\Delta}$ does not move when $S^{\Delta}$ steps in the direction $\leftrightarrow$, while it makes an increment -1 (resp., +1) when the increment of $S^{\Delta}$ is $\swarrow$ or $\searrow$ (resp., $\nearrow$ or $\nwarrow$). 
%Both skip some of the steps of $S^{\Delta}$: $S^{1,\Delta}$ when $S^{\Delta}$ steps in `$\nwsearrow$'' direction, and  when $S^{\Delta}$ steps in ``$\neswarrow$'' direction.  
Let $\tau^{i,\Delta}_{\pm L}$ be the hitting time of $\pm L$ by $S^{i,\Delta}$, and $M^i_n$ be the number of steps made by $S^{i,\Delta}$ by the time $S^{\Delta}$ makes $n$ steps.
%first times the corresponding projections walks hit either $+L$ or $-L$, and let $M^1_n$ and $M^2_n$ be the number of steps $S^{1,\Delta}$ and $S^{2,\Delta}$ make by the time $S^{\Delta}$ makes $n$ steps.
 Then we have
\begin{align}
 \P\left( \tau^{\Delta}_{L} > n \right) &= \P\left( \tau^{1,\Delta}_{\pm L/2} > n, \tau^{2,\Delta}_{\pm L/2} > n \right) \\
 &=\P\left( \tau^{1,\Delta}_{\pm L/2} > n, \tau^{2,\Delta}_{\pm L/2} > n\ |M^1_n \geq M^2_n \right)\times\P(M^1_n \geq M^2_n )\nonumber\\
 &\hspace{6cm}+ \P\left( \tau^{1,\Delta}_{\pm L/2} > n, \tau^{2,\Delta}_{\pm L/2} > n\ |M^1_n < M^2_n\right)\times\P(M^1_n < M^2_n)\nonumber\\
 &\leq \P\left( \tau^{1,\Delta}_{\pm L/2} > n\ |M^1_n \geq M^2_n \right)\times\P(M^1_n \geq M^2_n )+ \P\left(\tau^{2,\Delta}_{\pm L/2} > n\ |M^1_n < M^2_n\right)\times\P(M^1_n < M^2_n)\nonumber\\
 &\leq \P\left( \check{\tau}^{o}_{\pm L/2}>\frac{n}{2} \right).\nonumber
\end{align}
\Qed
\end{Lemma}

Theorem \ref{thm:ccp_lim} and the discussion preceding it suggest the following conjecture.
\begin{Conjecture}
\label{conjLine}
 The interface curve of the voter model in $B_{L}$ converges to a straight line as $L\rightarrow \infty$.
\end{Conjecture}  

Theorem \ref{thm:ccp_lim} also provides us with a useful test of the quality of our numerical estimation techniques (which are of course for finite $L$).  Having established that $\E[|C_L(o)|]\approx L^{\gamma}$ and $\E[|M_L|]\approx L^\beta$ with $\gamma=\beta=2$, we estimated the exponents from simulation data with estimators as in \eqref{dtilde},\eqref{dhat}, and \eqref{dregression} giving
\begin{align}
&\tilde{\gamma}=1.579,  \qquad \hat{\gamma}=1.861,  \qquad \gamma^*=1.841 \label{gamma}\\
&\tilde{\beta}= 1.657,   \qquad \hat{\beta}=1.911,   \qquad \beta^*=1.897. \label{beta}
\end{align}
Thus, again the $\hat{\bullet}$ estimators are closest to the true value.
%fitted the simple linear regression models
%\begin{align}
%\log(|M_L|)=&\beta\log(L)+\epsilon, \label{log_ML}\\
%\log(|C_L(o)|)=&\gamma\log(L)+\epsilon. \label{log_CL}
%\end{align}
%to simulated data with $L=128, 256, 512, 1024$.  
%With no apparent violation of the assumptions underlying such linear models, we obtained estimates %$\beta^*=1.897$ and $\gamma^*=1.841$ with standard errors $0.004$ and $0.01$ respectively (see Appendix ??). % {\bf If these are repeatable!!??  This is suggestive that the leading asymptotics of $L^2$ in terms of %expectation of these quantities has not kicked in at the size of $L$ we are looking at, or that the actual %sizes are not highly  concentrated around their mean???}

Figure \ref{fig:5largest} suggests that the (largest) coalescing classes are rather disconnected and sparse, which poses a potential problem for a rescaling argument like that mentioned at the beginning of Section \ref{sec:coalescing}. This is because the coalescing classes will not scale to single points if their diameters are $\ge cL$ with non-vanishing probability.  It is an open problem to prove that for some $c,c'\in (0,1)$, $\liminf_{L \ra \infty}\P(\exists x,y \in B_L: |x-y|>cL, \tau^{xy}<\tau_x\wedge\tau_y)\ge c'$.  This would imply that with positive probability there are coalescing classes with diameter at least $cL$.
%But it seems that the interface curve does not go around these clusters. 
A very large proportion of our simulated curves cut through $M_{L}$ in the sense that the interface curve has sites belonging to $M_{L}$ on both sides. The proportion increases from 0.9819 for the boxes of size 128,  %to 0.9946 for the boxes of size 256,  to 0.9973 for the boxes of size 512, and 
to 0.9985 for the boxes of size 1024.
%If we define connected clusters of a coalescing class to be the maximal connected subsets of coalescing %classes, 
Thus, the connected clusters/subsets $C^c_{L}(x)$ (containing $x\in B_L$) of coalescing classes $C_{L}(x)$,   may be better candidates to use in rescaling arguments as the interface curve has to go around them. 
%The comparison to the harmonic percolation breaks down in this case as the votes that connected clusters share are not independent any more because the votes of two connected clusters are the same if they are subsets of the same coalescing set. Yet the information on the connected clusters is easily extractable from the simulated samples and is interesting to look at on its own.
%For each configuration of a voter model in $B_{L}$ l
Assuming that  
$\E\left[ \left|C^c_{L}(x)\right|\right]\approx L^{\gamma'}$ and  $\E\left[ \max_{x\in B_L}\left|C^c_{L}(x)\right|\right]\approx L^{\beta'}$, we obtain the estimates $\hat{\gamma}'=1.548$ and $\hat{\beta}'=1.741$.  
%(although these appear to have not yet settled as functions of for our largest $L$).

%Under the same simplified assumption that $M^c_L\approx L^{\beta'}$ and $C^c_L(o)\approx L^{\beta''}$, using estimator for the exponents $\beta'$ and $\beta''$ analogous to the one used to estimate the exponent for expected length \eqref{dhat} and expected size of coalescing classes
%\begin{align*}
%  \hat{\beta}'_L = \log_2{\frac{\bar{M}^c_{2L}}{\bar{M}^c_{L}}} \text{ and } \hat{\beta}''_L = %\log_2{\frac{\bar{C}^c_{2L}(o)}{\bar{C}^c_{L}(o)}}, 
%\end{align*}
%we get the estimates that are substantially lower than 2 but are increasing with $L$ in a way that it is %hard to make any conclusions (see Appendix for more details). 

Another piece of information that may support Conjecture \ref{conjLine} is the displacement of the curve $Z_L$ from its conjectured diagonal limit, $D$. Assuming that 
%\begin{align*}
$\E\left[\max_{x\in Z_L}\min_{y \in D}|x-y|\right]\approx L^{\alpha}$,
%\end{align*}
%Assuming as usual that
%\begin{align*}
% D_L ,
%\end{align*}
%we estimate $\alpha$ using the usual 
we obtain the estimate $\hat{\alpha}=0.971$.
%\begin{align*}
% \hat{\alpha}_L = \log_2{\frac{\bar{D}_{2L}}{\bar{D}_{L}}},
%\end{align*}
%with $\alpha$ estimated around 0.97, which seems too close to 1, but still more than 10 standard errors away %from 1. For more detailed results see Appendix.

\section*{Acknowledgements}
MH thanks David Wilson for helpful discussions at the initial stages of this project.  MH and YM thank Raghu Varadhan, Federico Camia, and Pierre Nolin for helpful discussions and the Centre for eResearch at U.~Auckland for providing the computing resources and support.  The work of MH and YM was supported by an FRDF grant from U.~Auckland. The work of CMN was supported in part by US NSF grants OISE-0730136 and DMS-1007524. 

\bibliographystyle{plain}

\section{Appendix}
As indicated earlier, each of our estimates is based on 10000 independent simulations of the voter model for each value of $L$ being considered.  Although the simulations had a finite time horizon, in all cases all coalescing walks eventually reached the boundary.  All simulations were conducted in C++ and all statistical analyses and plots were performed in R.  The data is available on request, but at approximately 600GB, may be difficult to transfer.

Note also that all of our simulations actually took place on boxes of side length $L'=L+2$, so our estimators were actually
\begin{align}
\tilde{d}=\tilde{d}_{m,L}=&\frac{\log(\bar{Z}_m(L+2))}{\log(L)},\label{dtildeprime}\\
\hat{d}=\hat{d}_{m,L} =& \log_2\left(\frac{\bar{Z}_{m}(2L+2)}{\bar{Z}_m(L+2)}\right).\label{dhatprime}
%=\log_{2}\frac{Average\ length\ of\ curve\ in\ a\ box\ of\ size\ 2n}{Average\ length\ of\ curve\ in\ a\ box\ of\ size\ }.
\end{align}
Similarly our ordinary least squares estimator $d^*$ is in fact an estimator for the slope coefficient $d$ of the simple linear regression model
\begin{align}
\log(|Z_i|)=d\log(L_i)+a+\vep_i,\label{dregressionprime}
\end{align}
where $\{|Z_i|\}_{i\le n}$ are interface lengths on boxes of side lengths $\{L_i+2\}_{i \le n}$, and the $\vep_i$ are random variables with mean 0.  This does not change the consistency properties of the estimators, and e.g.~results in an estimate $\tilde{d}$ differing in only the fourth decimal place when we are dividing by $\log(L)$ (instead of $\log(L+2)$).

%\subsection{Interface length simulation}
%The results of simulations computing the lengths of the interface curves give the estimates $\hat{d},\tilde{d},d^*$ (with %$m=10000$) appearing in Table \ref{tab:d_estimates}. 

\subsection{The size of coalescing classes}
\begin{figure}
\includegraphics[scale=.5]{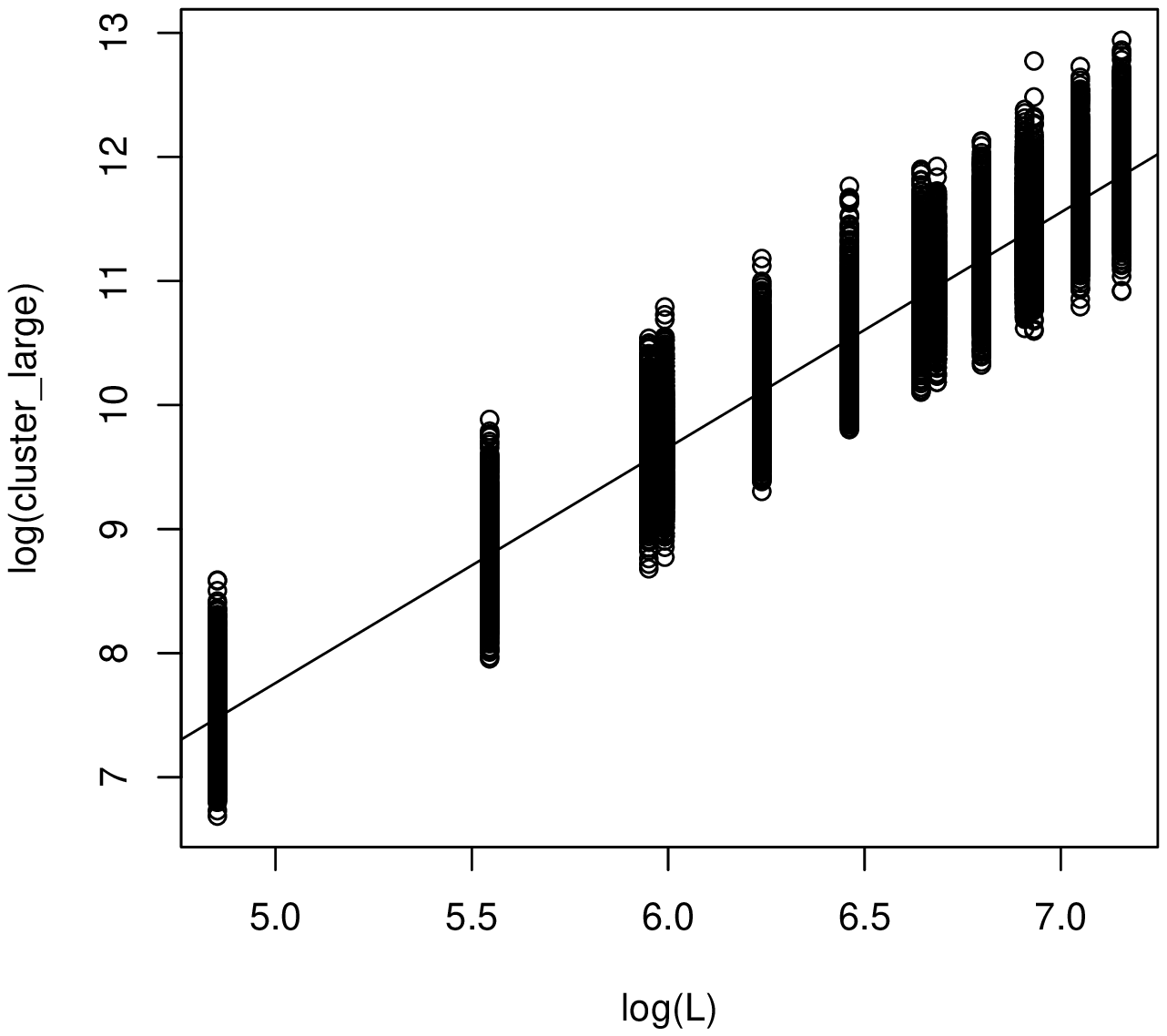}
\includegraphics[scale=.5]{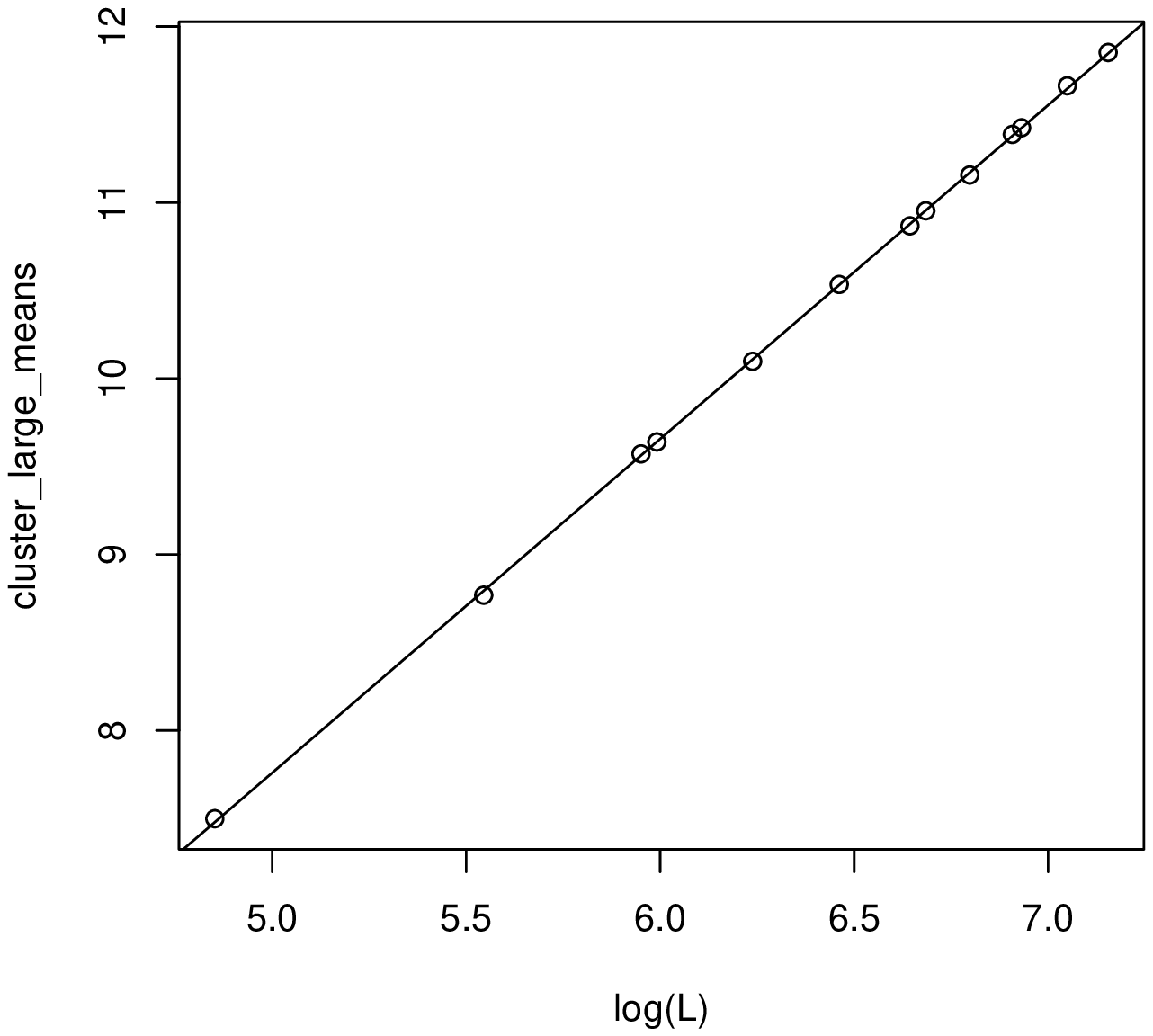}
\caption{The fitted line for \eqref{log_ML} for the largest coalescing class, with background the raw data and the fixed $L$ means respectively.}
\label{fig:ML}
\end{figure}

\begin{figure}
\includegraphics[scale=.5]{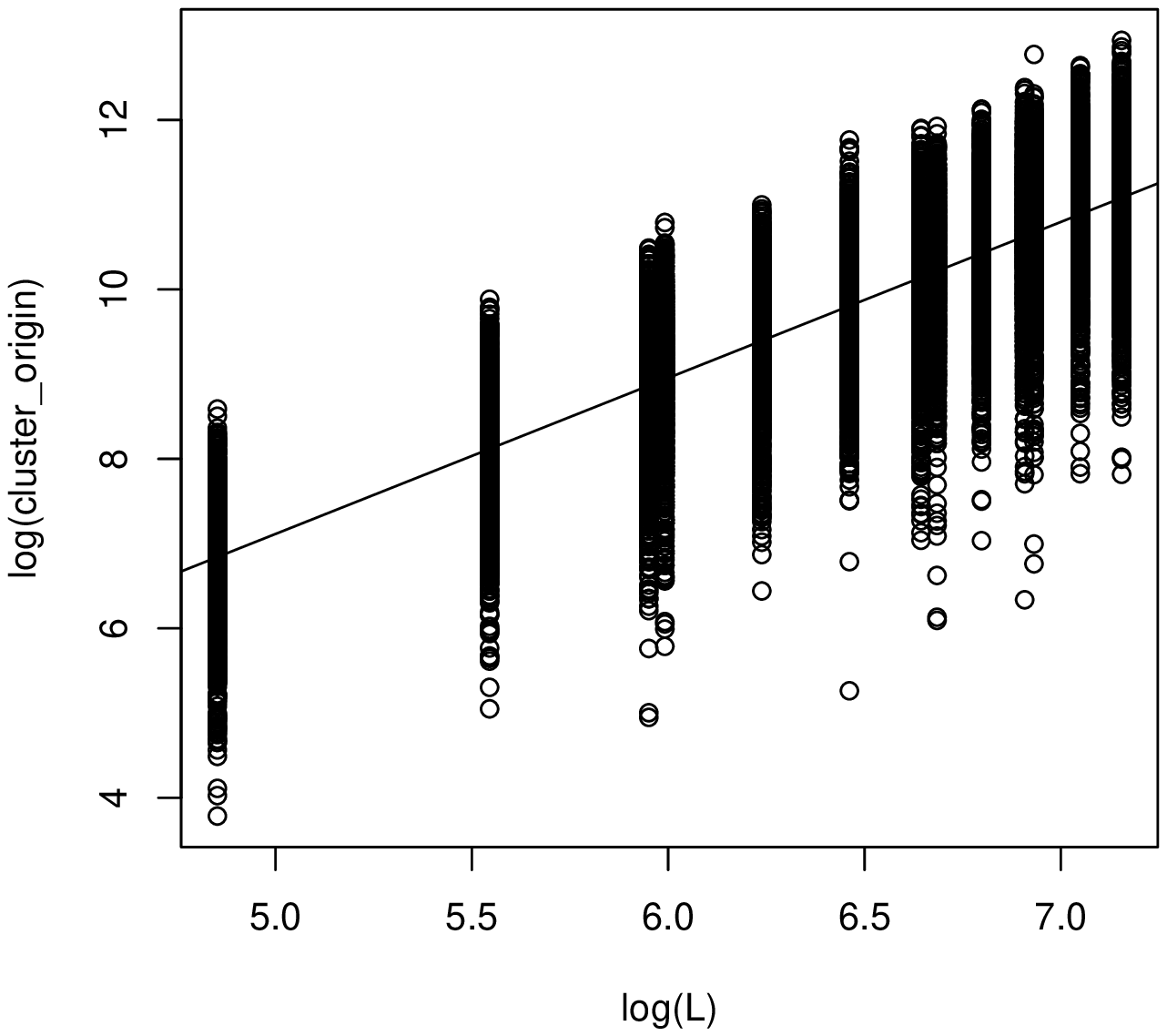}
\includegraphics[scale=.5]{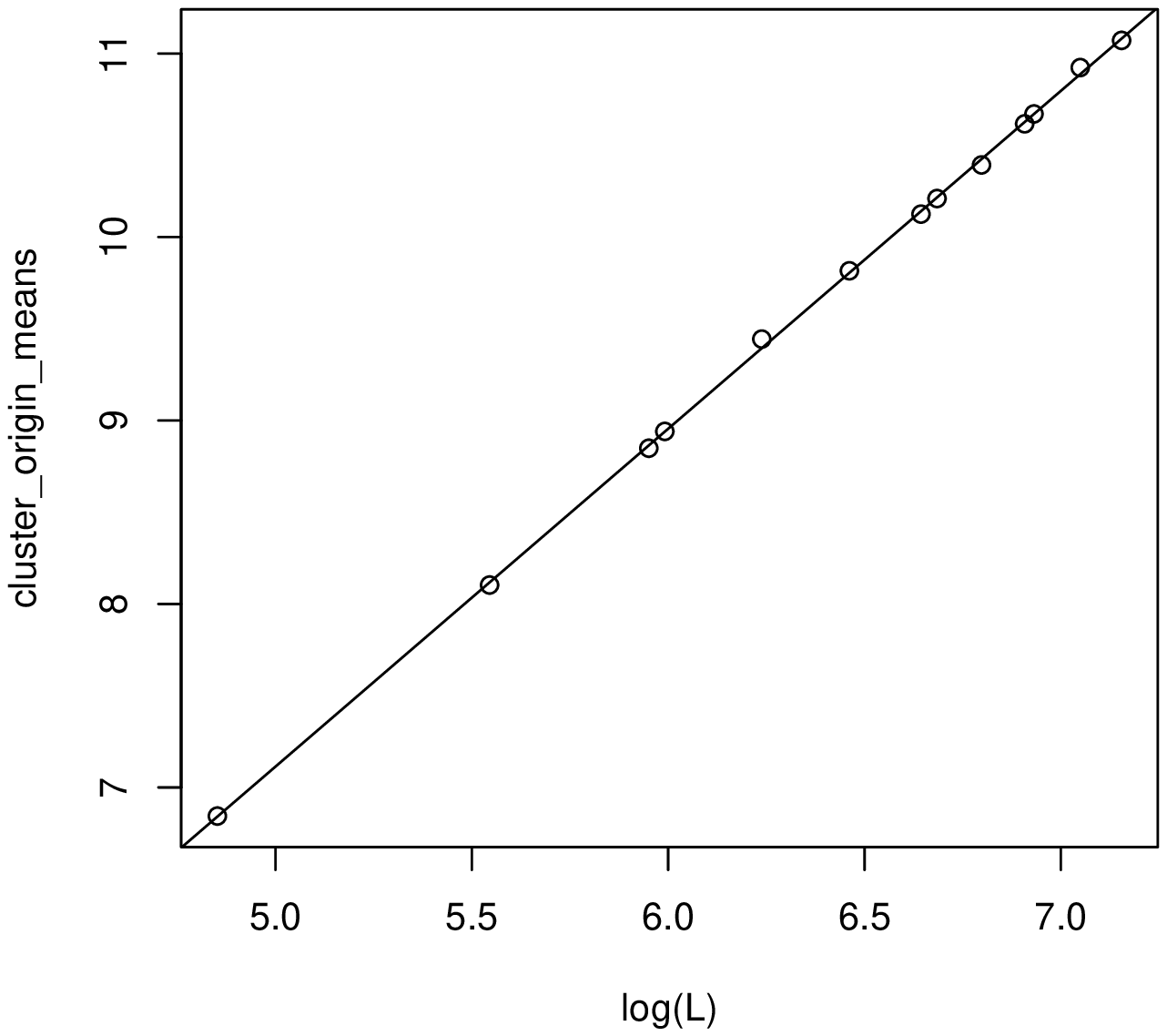}
\caption{The fitted line for \eqref{log_ML} for the coalescing class of a central vertex, with background the raw data and the fixed $L$ means respectively.}
\label{fig:origin_fit}
\end{figure}

Since the $\hat{\bullet}$ and $\tilde{\bullet}$ estimators are defined straightforwardly and have already been discussed somewhat at the end of Section \ref{sec:coalescing}, let us turn our attention here to the regression estimators $\bullet^*$ for the class sizes.  Assume that the assumptions prior to \eqref{dregression} hold for $|C_L(o)|$ and $|M_L|$ with exponents $\gamma$ and $\beta$ respectively, so that e.g. 
\begin{align}
\log(|M_L|_i)=\beta\log(L)+a+\vep_i. \label{log_ML}
\end{align}
We obtain an estimate $\beta^*$ of $\beta$ by fitting the simple linear model $\log(|M_L|)\sim\beta\log(L)$.  Fitting this linear model in R we obtain the following output:
\begin{verbatim}
Call:  lm(formula = log(cluster_large) ~ log(L))

Residuals:
     Min       1Q   Median       3Q      Max 
-0.93002 -0.22984 -0.02386  0.20920  1.34982 

Coefficients:
             Estimate Std. Error t value Pr(>|t|)    
(Intercept) -1.726828   0.028139  -61.37   <2e-16 ***
log(L)       1.897190   0.004374  433.72   <2e-16 ***

Residual standard error: 0.3201 on 12998 degrees of freedom
Multiple R-squared: 0.9354,     Adjusted R-squared: 0.9354 
F-statistic: 1.881e+05 on 1 and 12998 DF,  p-value: < 2.2e-16 
\end{verbatim}
Figure \ref{fig:ML} and the standard diagnostic tests suggests that the model fits very well.  However, the estimate for $\beta$ is more than 20 standard errors from the known (from Theorem \ref{thm:ccp_lim}) true value of $2$, so this estimator seems to be doing a poor job of estimating the true limiting behaviour in $L$.  We believe that this is caused by our inability to simulate the model for very large $L$.
%Using \eqref{log_ML} the simulations suggest that $M_L\approx L^{1.9}$, i.e., that $M_L$ grows slower than $L^{2}$. See also Figure \ref{fig:ML}.
For the coalescing class of the center, the linear regression estimator fits less well and gives an estimate of $\gamma^*=1.84$ (see also Figure \ref{fig:origin_fit}).
\begin{verbatim}
Call:  lm(formula = log(cluster_origin) ~ log(L))

Residuals:
    Min      1Q  Median      3Q     Max 
-4.5422 -0.4825  0.0914  0.5688  2.1029 

Coefficients:
            Estimate Std. Error t value Pr(>|t|)    
(Intercept) -2.09387    0.06964  -30.07   <2e-16 ***
log(L)       1.84149    0.01083  170.10   <2e-16 ***

Residual standard error: 0.7923 on 12998 degrees of freedom
Multiple R-squared:  0.69,      Adjusted R-squared:  0.69 
F-statistic: 2.894e+04 on 1 and 12998 DF,  p-value: < 2.2e-16 
\end{verbatim}

\subsection{Connected clusters of a coalescing class}
Recall from the end of Section \ref{sec:coalescing} that for each configuration of a voter model in $B_{L}$, $C^c_{L}(x)$ denotes the connected subset (containing $x$) of the coalescing class of $x$. Letting $M^c_L$ denote a largest such connected subset
and assume
\begin{align*}
\E\left[ \left|C^c_{L}(x)\right|\right]\approx L^{\gamma'} \text{ and } \E\left[| M^c_L|\right]\approx L^{\beta'}.
\end{align*}

Recalling the discussion around \eqref{dhat},
% and the fact that the $\hat{\bullet}$ estimators seem to be performing best, 
we obtain the following $\hat{\bullet}$ estimates, which indicate that this estimator converges more slowly with increasing $L$ than the corresponding estimator for the expected curve length exponent.

%\begin{table}[h!]
\begin{center}
\begin{tabular}{|l|l|l|l|}
\hline
$L$&128&256&512\\
\hline
$\hat{\beta'}_L$&1.69253&1.71810&1.74051\\
\hline
$\hat{\gamma'}_L$&1.47049&1.49151&1.54791\\
\hline
\end{tabular}
\end{center}
%\end{table}

\subsection{Maximum displacement of the curve from the diagonal}
Recall the last paragraph of Section \ref{sec:coalescing}.  Under the assumption that 
$\E\left[\max_{x\in Z_L}\min_{y \in D}|x-y|\right]\approx L^{\alpha}$ we have the following $\hat{\bullet}$ estimates for $\alpha$.
%\begin{align*}
% D_L \approx L^{\alpha},
%\end{align*}
%we estimate $\alpha$ by
%\begin{align*}
% \hat{\alpha}_L = \log_2{\frac{\bar{D}_{2L}}{\bar{D}_{L}}},
%\end{align*}
%with the following results
%\vspace{0.5cm}
%\begin{table}[h!]
\begin{center}
\begin{tabular}{|l|l|l|l|}
\hline
$L$ &128&256&512\\
\hline
$\hat{\alpha}_L$&0.96847&0.97545&0.97068\\
\hline
\end{tabular}
\end{center}
%\end{table}

%Here the value of about 0.97 for $\alpha$ is too close to 1 to reach any clear conclusion.
%But some statistical tests suggest that $\alpha < 1$. In the first of these we apply a linear regression analysis for the relation 
%\begin{align*}
% \log{D^i_L} = \alpha \log{L} +\epsilon^i, 
%\end{align*}
%where $D^i_L$ is an individual observation of maximum displacement whose mean is $D_L$, and $\epsilon^i$ is a random error. We %use 40000 points for all our curves of all sizes together. The results are the following

Fitting a simple linear model to the data gives a very similar estimate of 0.9689 (which is about 18 standard errors away from 1) as per the following.
\begin{verbatim}
Residuals:
     Min       1Q   Median       3Q      Max 
-1.11443 -0.17013  0.01359  0.18474  0.59812 

Coefficients:
             Estimate Std. Error t value Pr(>|t|)    
(Intercept) -1.090206   0.009573  -113.9   <2e-16 ***
log(Lvec)    0.968864   0.001611   601.4   <2e-16 ***

Residual standard error: 0.2497 on 39998 degrees of freedom
Multiple R-squared: 0.9004,     Adjusted R-squared: 0.9004 
F-statistic: 3.617e+05 on 1 and 39998 DF,  p-value: < 2.2e-16 
\end{verbatim}

\end{document}